%% file: main.tex
\documentclass[preprint,review,3p]{elsarticle}

\usepackage{amsmath}
\usepackage{amssymb}
\usepackage{braket,amsfonts}
\usepackage{array}
\usepackage{pgfplots}
\usepackage{algorithmic}
\usepackage{tikz}
\usepackage{graphicx,epstopdf}
\usepackage{placeins}
\usepackage[caption=false]{subfig}
\usepackage{lineno}

\journal{Journal of Computational Physics}

\input{abbrev.tex}
\begin{document}

\begin{frontmatter}

\title{A High Order Cartesian Grid, Finite Volume Method for Elliptic Interface Problems}

\author[inst1,inst2]{Will Thacher\corref{cor1}}
\cortext[cor1]{Corresponding author: wthacher@lbl.gov}

\affiliation[inst1]{organization={ Graduate Student Researcher, {Applied Science and Technology Group}, University of California Berkeley}, 
city={Berkeley},state={CA},
postcode={94720},
country={United States}}

\affiliation[inst2]{organization={ {Applied Numerical Algorithms Group}, Lawrence Berkeley National Laboratory}, 
city={Berkeley},state={CA},
postcode={94720},
country={United States} }

\author[inst2]{Hans Johansen}
\author[inst2]{Daniel Martin}

\begin{abstract}
We present a higher-order finite volume method for solving elliptic PDEs with jump conditions on interfaces embedded in a 2D Cartesian grid. Second, fourth, and sixth order accuracy is demonstrated on a variety of tests including problems with high-contrast and spatially varying coefficients, large discontinuities in the source term, and complex interface geometries. We include a generalized truncation error analysis based on cell-centered Taylor series expansions, which then define stencils in terms of local discrete solution data and geometric information. In the process, we develop a simple method based on Green's theorem for computing exact geometric moments directly from an implicit function definition of the embedded interface. This approach produces stencils with a simple bilinear representation, where spatially-varying coefficients and jump conditions can be easily included and finite volume conservation can be enforced. 
\end{abstract}


\end{frontmatter}



\makeatletter
\let\LN@align\align
\let\LN@endalign\endalign
\renewcommand{\align}{\linenomath\LN@align}
\renewcommand{\endalign}{\LN@endalign\endlinenomath}

\section{Introduction} \label{S:intro}
Elliptic PDEs with discontinuities in the source term, coefficients and solution form an important class of equations in computational science and engineering. These equations arise from mathematical models of multi-material systems, multi-phase flows, crystal growth, and many other physical processes \cite{LiReview}. Solving such equations numerically is not straightforward because the accuracy of the scheme is typically based on smoothness assumptions that do not in general apply at the interface. 

Numerous schemes have been proposed to solve this problem based on finite difference, finite volume, and finite element formulations. These methods can roughly be classified into those that treat the interface explicitly by creating elements that conform to the shape of the interface, or those that represent the interface implicitly by ``embedding" it onto a non-conforming mesh (see \cite{gibouReview} for a thorough review and further references). In the finite element realm, methods such as \cite{Babuska1970TheFE} body-fit the mesh to the interface whereas methods such as \cite{IIMFE} use a fixed mesh and modify basis functions where the interface crosses elements. A widely used and influential method in the finite difference category is the Immersed Interface Method (IIM) \cite{IIM}. The IIM uses standard Cartesian grid finite difference stencils away from the interface and modifies stencils near the interface using one-sided Taylor series expansions that incorporate jump conditions. The Ghost Fluid Method \cite{GFM} extrapolates the solution across the interface to nearby grid points by incorporating jump conditions so that standard stencils can still be used at all grid points. These various finite difference methods are closely related to various schemes for imposing boundary conditions; jump conditions can be thought of as a sort of implicit boundary condition that depends on the solution itself. 

This paper is concerned with the third category: finite volume schemes. These methods are conservative (in the sense that the divergence theorem is applied to a control volume), and have well-studied stability properties \cite{Hans_thesis}. The embedded boundary (EB) method of \cite{Crockett} combines the implicit and explicit interface representations: the interface is embedded onto a Cartesian grid, forming cut-cell volumes of arbitrary shapes where it intersects rectangular cells.  The elliptic equation is then discretized in flux divergence form using techniques developed in \cite{Hans_thesis} and \cite{SCHWARTZ2006531} with appropriate modifications made at the interface to enforce jump conditions. The method developed in \cite{Crockett} is second-order accurate $L^{\infty}$ and $L^1$ norm, but is difficult to extend to higher order accuracy. For low-order methods, cell averages can be treated as point values to second-order accuracy, so finite difference type schemes can be employed to create stencils. This is not the case for higher order finite volume methods; integration must be performed over arbitrarily-shaped ``cut cells.''  Techniques such as choosing a midpoint or centroid as a quadrature rule  for surface integrals of the flux is not sufficient for higher-order accuracy. 

Many of these difficulties are being addressed by recent developments in higher-order finite volume and EB methods, which are summarized in \cite{HOFV}. One example is the use of weighted least-squares interpolation for stencil construction in complex geometries (\cite{chen}, \cite{Dharshi}) as well as the derivation of high-order stencils on Cartesian grids \cite{Zhang}. Given the close relationship between boundary conditions and jump conditions, we propose extending the methodology of \cite{Dharshi} to the elliptic interface problem. The primary contribution of this research is a finite volume method for the variable coefficient 2D elliptic interface problem that is 1) high order accurate and 2) conservative. In the process, we have also created approaches for 1) an efficient technique for generating exact geometric information from an implicit function and 2) a method for building high-order finite volume stencils for variable coefficient elliptic operators on arbitrary cut-cell meshes. 

The outline of the paper is as follows: In Section 2, we define mesh and geometric quantities and give a general truncation error analysis which allows us to design stencils of arbitrarily high order. In Section 3, we describe in detail our method for constructing stencils. In Section 4, we present results that validate the approach using a series of model problems that test different aspects of the scheme.

\section{Discretization}\label{sec:GE}
\begin{figure}
    \centering
    \begin{tikzpicture}
    \filldraw[color=black,fill=yellow!50,thick] (0,0) rectangle (5,5);
    \filldraw[color=red,fill=green!40,very thick] (5,2.5) arc[start angle=270, end angle=180, radius=2.5];
    \fill[color=green!40] (5,2.5) -- (5,5) -- (2.5,5) -- cycle;
    \draw [-stealth, very thick, black](3.2322330470336307,3.232233047033631) -- (3.75,3.75);
    \draw[color=black,thick](0,0) rectangle (5,5);
    
    \node[] at (1.5,1.5) {$V_{+, \*i}$};
    \node[] at (4,4) {$V_{-, \*i}$};
    \node[] at (2.2,4) {$A_{B, \*i}$};
    \node[] at (3.05,3.2) {$\hat{\*n}$};
    \node[] at (2.5,-.5) {$A_{+,\*i - \frac{1}{2} \*e_y}$};
    \node[] at (1.25,5.5) {$A_{+,\*i + \frac{1}{2} \*e_y}$};
    \node[] at (3.75,5.5) {$A_{-,\*i + \frac{1}{2} \*e_y}$};
    \node[] at (-1,2.5) {$A_{+,\*i - \frac{1}{2} \*e_x}$};
    \node[] at (6,1.25) {$A_{+,\*i + \frac{1}{2} \*e_x}$};
    \node[] at (6,3.75) {$A_{-,\*i + \frac{1}{2} \*e_x}$};
    \end{tikzpicture}
    \caption{Cut cell geometric quantities that make up the finite volume notation.}
    \label{fig:cutCell}
\end{figure}
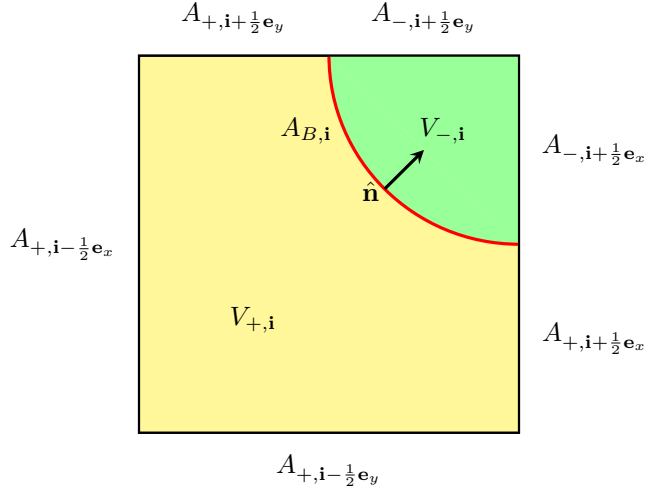

Let $\Omega$ be a physical domain that is divided into the subdomains $\Omega^+$ and $\Omega^-$ by an interface $\Gamma$. We consider the variable coefficient elliptic interface problem for $u(x)$:
\begin{align}
    \alpha u - \nabla \cdot \lrp{\beta \nabla u} &= f \label{eq:elliptic} \hbox{ on } \Omega \\
    \lrb{u} &= w \hbox{ on } \Gamma \label{eq:uJump}\\
    \lrb{\beta \partial_{\*n} u} &= v \hbox{ on } \Gamma \label{eq:dujump} \, . 
\end{align}
Here, $\lrb{\cdot}$ denotes a jump in some quantity at the interface: $\lrb{u} = u^+(\*x) - u^-(\*x)$ at some point $\*x \in \Gamma$, and the term $\partial_{\*n} u \equiv \nabla u \cdot \hat{\*n}$ represents a \textit{flux} at this boundary with unit normal $\hat{\*n}(\*x)$. Finally, coefficients $\alpha^\pm(\*x), \beta^\pm(\*x)$ and the source term $f^\pm(\*x)$ vary in space and may be discontinuous across $\Gamma$.

The domain $\Omega$ is discretized into a Cartesian mesh of square control volumes (or  ``cells") $V_{p,\*{i}}$, $\*i \in \mathbb{Z}^2$, that have centroids $\*x_{p,\*i}$ and side lengths of scale $h$, the grid spacing (see Figure \ref{fig:cutCell}). We indicate $p \in \set{+,-}$ to specify a subdomain of $\Omega_p$; $p$ can often be thought of as the phase or material type of a physical quantity. We assume that each cell $V_{p,\*i}$ may have up to four grid-aligned faces, which we label $A_{p,\*i \pm \frac{1}{2} \*e_d}$, where $\*e_d$ is the unit vector in direction $d$. 

Any cell that is intersected by $\Gamma$, the ``embedded boundary'' (EB), is called a ``cut" cell. We make the following assumptions to simplify the geometric considerations. First, a cut cell consists of only two control volumes $V_{+,\*i}$ and $V_{-,\*i}$ divided by a portion of the EB, denoted by $A_{B,\*i}$. The unit normal vector $\hat{\*n}$ on $\Gamma$ points from $\Omega^+$ to $\Omega^-$. So, along with $A_{p,\*i \pm \frac{1}{2} \*e_d}$ as the grid-aligned faces of each portion of the cut cell, each cut cell must have a total of at least 3, and at most 5, faces.

Because we are using a finite volume formulation, we should define additional geometric quantities that will be useful throughout this paper: a geometric ``moment" is an integral of a centered monomial over some specified region. We define four moments corresponding to four components of the geometry:
\begin{align}
    m_{{p,\*i}}^{\*q} &= \int_{V_{p,\*i}} \*{(x-\bar{x})^q} \ dV \label{eq:vMom}\\
    m_{{p,\*i \pm \frac{1}{2}e_d}}^{\*q} &= \int_{A_{p,\*i \pm \frac{1}{2}e_d}} \*{(x-\bar{x})^q} \ dA  \label{eq:aMom}\\
    m_{{B,\*i}}^{\*q} &= \int_{A_{B,\*i}} \*{(x-\bar{x})^q} \ dA \label{eq:EBMom}\\
    m_{{B,\*i,d}}^{\*q} &= \int_{A_{B,\*i,d}} \*{(x-\bar{x})^q} \ \hat{n}_d \ dA  \ , \label{eq:nMom}
\end{align}
where $\*q =\lrb{q_x \ q_y}$ is a vector of non-negative integers, and we use the multi-index notation $\*{(x-\bar{x})^q} = (x-\bar{x})^{q_x} (y-\bar{y})^{q_y}$. The multi-indices have sum of at most $P$, $|\*q| \le P$, and are ordered lexicographically: $\set{00,10,20,...P0,01,11, ... 0P}$. This allows us to refer to $\*v \lrb{\*q}$ as the $\*q^{th}$ entry of a vector $\*v$. Thus by definition, the volume of cell $\*i$ is $\abs{V_{p,\*i}} = m_{{p,\*i}}^{00}$, and the centroid $\bar{\*x}_{p,\*i}$ of $V_{p,\*i}$ is $\frac{1}{|\Vpi|}\lrb{m_{{p,\*i}}^{10}, m_{{p,\*i}}^{01}}$. Similarly, $m_{{B,\*i}}^{00}$ is the area of the EB, and $m_{{B,\*i,x}}^{00}$ is its $x$ normal component-weighted area, or $x$ direction \emph{cross-section}. For ease of notation, throughout this paper we ignore $\bar{\*x}$, although in practice it is the cell-center of each full Cartesian cell. Finally, our method for calculating these moments in 2D is exact to within roundoff errors, and is detailed in the appendix.

Two types of variables are stored on the mesh: cell-averaged quantities $\IP{ u }_{p,\*i} = \frac{1}{|\Vpi| } \int_{\Vpi} u \ dV $, and centroid-centered quantities $u_{p,\*i} = u(\*x_{p,\*i})$. The coefficients $\alpha$ and $\beta$ are given as point values at the centroids of cells. The right-hand-side function is provided as cell-averaged values, $\avg{f}$, of sufficient accuracy, and we solve for cell-averaged values $\avg{u}$.

Within this context, our finite volume scheme solves the discrete system:

\begin{align}
    \IP{\alpha  u}_{p,\*i} - \IP{ \nabla \cdot \beta \nabla u}_{p,\*i} &= \IP{f}_{p,\*i} \label{eq:ellipticD}\\
    \lrb{u}_{\*i} &= w_{\*i} \label{eq:uJumpD} \\
    \lrb{\beta \partial_{\*n} u}_{\*i} &= v_{\*i} \ . \label{eq:duJumpD}
\end{align}
for $\IP{u}_{p,\*i}$ in each volume $\Vpi$ in the mesh, subject to problem-specific boundary conditions. The system of equations we solve will have one degree of freedom in full cells and two degrees of freedom in cut cells. $\lrb{\cdot}_{\*i}$ denotes the integral of the jump of a quantity across the EB in cut cell $\*i$:
\begin{align}
    \lrb{u}_{\*i} = \int_{A_{B,\*i}} u^+ - u^- dA \ .
\end{align}

The objective of the following section is to provide a general truncation error analysis that will allow us to discretize \eqref{eq:ellipticD} -- \eqref{eq:duJumpD} to high order accuracy.

\subsection{Error Analysis}
In the finite volume or finite difference context, a \emph{stencil} approximates some functional $G(u)$ by a linear combination of local information about the function $u$. This functional is typically point values or integrals over some region of derivative of $u$. The function information, or data, can include boundary conditions, jump conditions, point values of the function, or averages of the function over some nearby region. Let $\*d$ denote this vector of local function data, and $\*s$ be the vector of stencil values corresponding to each of the data points in $\*d$. The truncation error $\tau$ of this stencil is defined as:
\begin{align}
    \tau = \*s^T \*d - G(u) \ .\label{eq:tau}
\end{align}
For the present problem, $G$ will be an integral of $u$, or some combination of its partial derivatives, over a one or two dimensional region. If we can approximate $u$ using a truncated Taylor series, $G(u)$ can be approximated with a linear combination of Taylor series coefficients of $u$ up to the desired order of accuracy. 

Throughout this section we drop the subscripts $p$ and $\*i$ except where necessary for clarity. We can express $u$ locally as a Taylor series expansion and remainder term:
\begin{align}
    u(\*x) = \sum_{|\*q| \leq P} \frac{1}{\*{q!}} u^{(\*q)} (\*{\bar{x}}) \*x^{\*q} + O\lrp{h^{P+1}} , \label{eq:taylor}
\end{align}
where, again using multi-index notation, $\*x^{\*q}= x^{q_x} y^{q_y}$, $\*{q!} = q_x! q_y!$, and $u^{(\*q)} = \frac{\partial^{q_x} \partial^{q_y}  }{\partial x^{q_x} \partial y^{q_y}} u $. The Taylor polynomial is then just
\begin{align}
    u(\*x) & = \sum_{|\*q| \leq P} c^{\*q}_{u} \*x^{\*q} + O\lrp{h^{P+1}} \label{eq:taylorcoef} \hbox{, where } \\
    c^{\*q}_{u} & = \frac{1}{\*{q!}} u^{(\*q)} (\*{\bar{x}}) \, . \nonumber
\end{align}

Integrating the flux divergence term in \eqref{eq:ellipticD} over the discrete volume $V$ and applying Gauss' theorem we obtain:
\begin{align}
    \int_{V} \nabla \cdot \beta \nabla u \ dV &= \lrb{ \sum_{\pm, \ d} (\pm 1) \int_{A_{\pm \frac{1}{2}e_d}} \beta \frac{\partial u}{\partial x_d} \ dA } + \int_{A_{B}} \beta \nabla u \cdot \hat{\*n} \ dA \ .\label{eq:fluxDiv}
\end{align}

Expressing $\beta$ and $u$ as Taylor expansions, for a surface integral over any (EB or grid-aligned) face $A$ we have:
\begin{align}
\label{eq:fluxint}
    \int_{A} \beta \nabla u \cdot \hat{\*n} \ dA &= \int_A \left( \sum_{|\*r| \leq P}  c_{\beta}^{\*r} \*x^{\*r} \right) \left( \sum_{|\*q| \leq P} c_{u}^{\*q} \lrb{ \frac{\partial \*x^{\*q}}{\partial x} \hat{n}_x + \frac{\partial \*x^{\*q}}{\partial y} \hat{n}_y} \right) + O(h^P) \ dA \\
    &= \int_A \left( \sum_{|\*r| \leq P}  c_{\beta}^{\*r} \*x^{\*r}  \right) \left( \sum_{|\*q| \leq P} c_{u}^{\*q} \lrb{ q_x\*x^{\*q-\*e_x} \hat{n}_x + q_y\*x^{\*q-\*e_y} \hat{n}_y} \right) + O(h^P) \ dA \\
    &= \sum_{|\*q|  \leq P}  \lrb{\sum_{|\*r|  \leq (P - |\*q|)}  c_{\beta}^{\*r} \lrp{ q_x m_{A,x}^{\*q + \*r -\*e_x} + q_y m_{A,y}^{\*q + \*r -\*e_y} }} c_{u}^{\*q} + O(h^{P+1}) \ .\label{eq:fluxTaylor}
\end{align} 
In order to compute $\IP{ \nabla \cdot \beta \nabla u} $, we apply \eqref{eq:fluxTaylor} to each surface integral in \eqref{eq:fluxDiv} and divide by the cell volume. Since $\abs{V}$ is an $O(h^2)$ quantity, the error term for the cell-averaged flux divergence term is $ O(h^{P-1})$. \\

Similarly, for the linear term in \eqref{eq:ellipticD}, we have:
\begin{align}
    \IP{\alpha  u} &=
    \frac{1}{\abs{V}} \int_{V} \sum_{|\*r| \leq P}  c_{\alpha}^{\*r} \*x^{\*r} \sum_{|\*q| \leq P}  c_{u}^{\*q} \*x^{\*q}  + O(h^{P+1}) dV \\
    &= \frac{1}{\abs{V}} \sum_{|\*q| \leq P} \lrb{\sum_{|\*r| \leq (P - |\*q|-2)} c_{\alpha}^{\*r} m^{\*r + \*q} } c_{u}^{\*q}  +O(h^{P+1})  \ .\label{eq:alphaTS}
\end{align}

 We see that our approximation to the functionals $G$ of interest can be written in the general form:
\begin{align}
    G(u) = \sum_{\abs{\*q} \leq P}  g^{\*q} c_{u}^{\*q} + O(h^R) = \*g^T \*c_{u} + O(h^R) \ , \label{eq:gtc}
\end{align}
where $\*c_{u}$ is the vector of approximate Taylor series coefficients for $u$, and $\*g$ is the vector of terms that result from operating on these coefficients, shown in brackets in \eqref{eq:fluxTaylor} and \eqref{eq:alphaTS}. The order $R$ error term for the functional approximation may be different than the order $P+1$ error term for the Taylor series of $u$ as a result of differentiation and integration.

To calculate the truncation error $\tau$ in \eqref{eq:tau}, we must fill in the vector $\*d$ with function data. Suppose, for example, that the function data are cell-averaged values of $u$. Then we can write:
\begin{align}
    \IP{u}_{\*j} &= \frac{1}{|V_{\*j}|} \int_{V_{\*j}} \sum_{|\*q| \leq P}  c^{\*q}_{u}  \*x^{\*q} + O\lrp{h^{P+1}} dV \\
    &= \sum_{|\*q| \leq P}  \frac{m_{\*j}^{\*q}}{|\Vj|}    c^{\*q}_{u} + O\lrp{h^{P+1}} = \*m_{\*j}^T \*c_{u} + O\lrp{h^{P+1}} \ ,
\end{align}
where $ \*m_{\*j}$ is the vector of cell-averaged volume moments for any cell $\*j$. If we have $n$ such cell-averaged values, then we can write:
\begin{align}
    \*d = \begin{bmatrix}
    \*m_{\*j_1}^T \\
    \*m_{\*j_2}^T \\
    . . . \\
    \*m_{\*j_n}^T
    \end{bmatrix} \*c_{u} + O(h^{P+1}) = \*M \*c_{u} + O(h^{P+1}) \ , \label{eq:MPhi}
\end{align}
where $\*M$ is the ``moment matrix" whose rows are the vectors $\*m_{\*j_i}^T$. Inserting this into \eqref{eq:tau} we have:
\begin{align}
    \tau = \*s^T \*M \*c_{u} + O(h^{P+1}) - \lrp{ \*g^T \*c_{u} + O(h^R) } \ .
\end{align}
To ensure that $\tau$ is $O(h^{\min(R,P+1)})$, we must have:
\begin{align}
    \*s^T \*M \*c_{u} &= \*g^T \*c_{u} \, .
\end{align}
This must hold for any $u$ with a Taylor series represented by arbitrary $\*c_{u}$, so that
\begin{align}
    \*M^T \*s &= \*g \, . 
\end{align}
Note that this linear system is a general form for the \emph{method of undetermined coefficients}, as explained in \cite{leveque} (Chapter 1.2): stencil weights are chosen so that the sum of Taylor series terms of the stencil exactly matches the Taylor series terms of the functional up to some order. In this case, we are working with cell averages rather than pointwise evaluations, but the principle is the same. If the neighborhood of local function data is chosen such that this linear system is exactly determined or undetermined, the least norm solution to this system is given by the pseudoinverse of $\*M$:
\begin{align}
    \*s = (\*M^T)^{\dagger} \*g \, . \label{eq:stencilForm}
\end{align}
If this system is undetermined, there are infinitely many stencils that will have the same order of truncation error. However, there is another profitable way to view this stencil construction process that justifies using the least-norm stencil. Suppose we want to interpolate the data stored in the vector $\*d$ with a degree $P$ polynomial. If the matrix $\*M$ is full rank, we can determine the coefficients of this polynomial by solving, in a least-squares sense, the linear system in \eqref{eq:MPhi}:
\begin{align}
    \*c_{u} &= \*M^{\dagger} \*d  \ . \label{eq:taylorLS}
\end{align}
Inserting this into \eqref{eq:gtc}, we have:
\begin{align}
    G(u) &= \*g^T \*M^{\dagger} \*d + O\lrp{h^{\min (R,P+1)}}  \ .
\end{align}
Finally, inserting this into the truncation error expression \eqref{eq:tau}:
\begin{align}
    \*s^T \*d &= \*g^T \*M^{\dagger} \*d \implies \\
    \*s &= (\*M^T)^{\dagger} \*g \ , \label{eq:sten}
\end{align}
where the implication follows from the fact this must hold for any $\*d$. 

Thus we can view our stencils as originating from either the undetermined system $\*M^T \*s = \*g$, in which $\*s$ is chosen to cancel lower order Taylor series terms, or the overdetermined system $\*M \*c_{u} = \*d$, in which we fit an interpolating polynomial to local data. This error analysis is quite general; we have made no mention of the shape of the volumes in the mesh, only that we know their geometric moments to sufficient accuracy. The interface jump conditions and boundary conditions are considered to be pieces of function data that can be used to build stencils near an interface or boundary. Taking the overdetermined perspective, by doing so we will constrain our interpolating polynomials, and therefore the numerical solution, to match boundary and jump conditions.



\section{Stencil Construction}
In the previous section we showed that if we approximate the solution $u$ with a polynomial whose coefficients are mapped from local function data by $\*M^{\dagger}$, then our stencil will take the simple form \eqref{eq:stencilForm}. In this section we will describe in detail our method for computing the moment matrix $\*M$ and the vector $\*g$ of terms that result from operating on the Taylor polynomials approximating $u$. From \eqref{eq:fluxTaylor} and \eqref{eq:alphaTS}, we see that to achieve a truncation error of order $P-1$ for the two terms in \eqref{eq:ellipticD}, we need to ``calculate" Taylor coefficients $c_{\phi}^{\*q}$ for $\phi \in \set{u,\beta, \alpha}$ up to order $P$. In this section we again drop the subscript $p$ except where necessary for clarity.


Both the linear term and the flux divergence term can be expressed as linear combinations of the Taylor series terms $c_{u}^{\*q}$, as seen in \eqref{eq:alphaTS} and \eqref{eq:fluxTaylor}. For the linear term, let $\*g_{\alpha, \*i}$ be the vector whose $\*q^{th}$ entry is 
\begin{align}
   \*g_{\alpha, \*i}\lrb{ \*q} =  \sum_{|\*r| \leq (P - |\*q|-2)} c_{\alpha}^{\*r} m_{\*i}^{\*r + \*q}  \ .
\end{align}
Applying \eqref{eq:sten}, we have that a stencil for the linear term is given by:
\begin{align}
    \*s_{\alpha,\*i} =  (\*M_{u, \*i}^T)^{\dagger} \*g_{\alpha, \*i} \label{eq:stencilAlpha}
\end{align}
where the subscript $\alpha, \*i$ indicates a stencil for the linear term, $\IP{\alpha u}$, over cell $\*i$. For the integral of $\beta \nabla u \cdot \hat{\*n}$ over a face $A$, we let $\*g_{\beta, A}$ be the vector whose $\*q^{th}$ entry is
\begin{align}
    \*g_{\beta, A} \lrb{\*q} = \sum_{|\*r| \leq (P - |\*q|)}  c_{\beta}^{\*r} \lrp{ q_x m_{A,x}^{\*q + \*r -\*e_x} + q_y m_{A,y}^{\*q + \*r -\*e_y} } \ .
\end{align}
A stencil for the flux surface integral \eqref{eq:fluxint} is then given by: 
\begin{align}
  \*s_{\beta, A} = (\*M_{u, \*i}^T)^{\dagger} \*g_{\beta, A}  \ . \label{eq:stencilbeta} 
\end{align}
 Each entry in $\*g_{\alpha, \*i}$ involves the Taylor coefficients of the $\alpha$ field, and likewise each $\*g_{\beta, A}$ entry of involves the Taylor coefficients of the $\beta$ field. We must generate these coefficients for each volume $V$, which can be done with the same moment matrix formulation that we use to compute the Taylor coefficients $\*c_u$.

Recall that $\alpha$ is given as point values at the centroids of cells. Let $\mathcal{N}_{p,\*i}$ be some neighborhood of cells $\*j_k \in \left\{\*j_1, \dots ,\*j_n\right\}$ in phase $p$ around volume $\Vpi$. A cell $V_{p,\*j}$ in this neighborhood has centroid $\*x_{p,\*j}$, with $\alpha_{p,\*j} = \alpha(\*x_{p,\*j })$. Our function data  will consist of $\alpha_{p,\*j}$ at each point $\*j_k$ in the neighborhood, which we compile into the vector $\*d_{\alpha, p, \*i}$. Let $\*M_{\alpha,p, \*i}$ be the variable coefficient interpolation matrix whose rows consist of monomials in the Taylor expansion of $\alpha$ evaluated at $\*x_{p,\*j }$:
\begin{align}
    \*M_{\alpha,p, \*i} \lrb{k, \*q} = \*x_{p,\*j_{k} }^{\*q} \ .
\end{align}
Letting $\*c_{\alpha,p}$ be the vector of Taylor series coefficients, we have 
\begin{align}
    \*c_{\alpha,p} = \*M_{\alpha,p, \*i}^{\dagger} \*d_{\alpha, p, \*i} \ .
\end{align}
We can do the same for $\beta$ to calculate $\*c_{\beta,p}$.

Therefore we have that:
\begin{align}
    \*g_{\alpha, \*i} &= \*G_{\alpha, \*i}  \*M_{\alpha, \*i}^{\dagger} \*d_{\alpha, \*i} \ ,
\end{align}
where $\*G_{\alpha, \*i} $ is the matrix whose $\*q, \*r$ entry is $m_{\*i}^{\*r + \*q}$. Combining this with \eqref{eq:stencilAlpha}, we can finally write:
\begin{align}
    \*s_{\alpha, \*i}^T &= \*d_{\alpha, \*i}^T (\*M_{\alpha, \*i}^{\dagger})^T \*G_{\alpha, \*i}^T \*M_{u, \*i}^{\dagger} \ .
\end{align}

This stencil is bilinear in $\alpha$ and $u$, which is appropriate for the bilinear functional we are trying to approximate. We can obtain a similar formulation for the flux stencil $\*s_{\beta, A}^T$ for each face $A$ of the cell. Thus constructing stencils is just a matter of constructing the moment matrices for $u,\beta$ and $\alpha$. We now describe how to construct moment matrices at, near, and away from the interface.

\begin{figure}
    \centering
    \begin{tikzpicture}
    \fill[white!50] (0,0) rectangle (5,5);
    \fill[blue!100] (0,3) rectangle (1,4);
    \fill[blue!100] (1,3) rectangle (2,4);
    \fill[blue!100] (2,3) rectangle (3,4);
    \fill[blue!100] (2,2) rectangle (3,3);
    \fill[blue!100] (3,2) rectangle (4,3);
    \fill[blue!100] (3,1) rectangle (4,2);
    \fill[blue!100] (4,1) rectangle (5,2);
    \fill[blue!100] (4,0) rectangle (5,1);

    \fill[blue!50] (0,4) rectangle (1,5);
    \fill[blue!50] (1,4) rectangle (2,5);
    \fill[blue!50] (2,4) rectangle (3,5);
    \fill[blue!50] (3,3) rectangle (4,4);
    
    \fill[blue!50] (4,2) rectangle (5,3);
    \fill[blue!50] (0,2) rectangle (1,3);
    \fill[blue!50] (1,2) rectangle (2,3);
    \fill[blue!50] (2,1) rectangle (3,2);
    \fill[blue!50] (3,0) rectangle (4,1);
    \draw[step=1.0,black,thin] (0,0) grid (5,5);
    \draw[red,thick] (4.5,0) .. controls (3,3.8) and (1.5,3.6) .. (0,3.5);
    \end{tikzpicture}
    \caption{For the $P=2$ scheme, the regular cell footprint is a standard five-point Laplacian, and if any point in the footprint is a cut cell, it is then ``irregular.''  Cut cells are shown with dark shading, irregular cells with light shading, and the remaining white cells are ``regular.''}
    \label{fig:cellTag}
\end{figure}
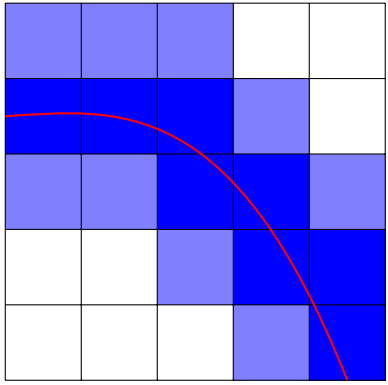

\subsection{Moment Matrices}
 We partition our cells into three subsets: cut cells $\Omega_C$, irregular cells $\Omega_I$, and regular cells $\Omega_R$. Cut cells are intersected by the EB. Irregular cells are not intersected by the EB, but at least one cell in the stencil footprint for a regular cell is intersected by the EB. See Figure \ref{fig:cellTag}. The method for constructing the moment matrices is different for each of these three types of cells.

\subsubsection{Regular Cells} \label{sec:regCellSten}

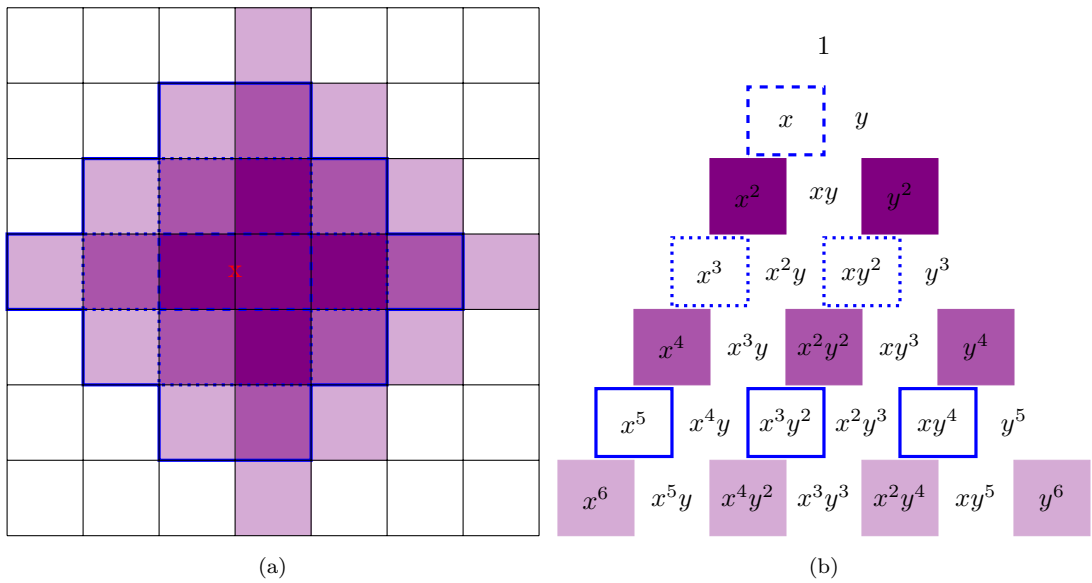
\begin{figure}
    \centering
    \subfloat[]{
\begin{tikzpicture}
    \fill[violet!100] (3,3) rectangle (4,4);
    \fill[violet!100] (2,3) rectangle (3,4);
    \fill[violet!100] (4,3) rectangle (5,4);
    \fill[violet!100] (3,4) rectangle (4,5);
    \fill[violet!100] (3,2) rectangle (4,3);

    \fill[violet!67] (1,3) rectangle (2,4);
    \fill[violet!67] (5,3) rectangle (6,4);
    \fill[violet!67] (3,1) rectangle (4,2);
    \fill[violet!67] (3,5) rectangle (4,6);
    \fill[violet!67] (2,4) rectangle (3,5);
    \fill[violet!67] (4,4) rectangle (5,5);
    \fill[violet!67] (2,2) rectangle (3,3);
    \fill[violet!67] (4,2) rectangle (5,3);
    
    \fill[violet!33] (0,3) rectangle (1,4);
    \fill[violet!33] (6,3) rectangle (7,4);
    \fill[violet!33] (3,0) rectangle (4,1);
    \fill[violet!33] (3,6) rectangle (4,7);
    \fill[violet!33] (2,5) rectangle (3,6);
    \fill[violet!33] (4,5) rectangle (5,6);
    \fill[violet!33] (2,1) rectangle (3,2);
    \fill[violet!33] (4,1) rectangle (5,2);
    \fill[violet!33] (5,2) rectangle (6,3);
    \fill[violet!33] (5,4) rectangle (6,5);
    \fill[violet!33] (1,2) rectangle (2,3);
    \fill[violet!33] (1,4) rectangle (2,5);

    \draw[very thick, dashed, blue] (2,3) -- (2,4) -- (4,4) -- (4,3) -- cycle;

    \draw[very thick, dotted, blue] (1,3) -- (1,4) -- (2,4) -- (2,5) -- (4,5) -- (4,4) -- (5,4) -- (5,3) -- (4,3) -- (4,2) -- (2,2) -- (2,3) -- cycle;

    \draw[very thick, blue] (0,3) -- (0,4) -- (1,4) -- (1,5) -- (2,5) -- (2,6) -- (4,6) -- (4,5) -- (5,5) --  (5,4) -- (6,4) -- (6,3) -- (5,3) -- (5,2) -- (4,2) -- (4,1) -- (2,1) -- (2,2) -- (1,2) -- (1,3) -- cycle;
    
    \node[red, very thick] at (3,3.5) {x};
    
    \draw[step=1.0,black,thin] (0,0) grid (7,7);
    
\end{tikzpicture}  
}
\subfloat[]{
    \begin{tikzpicture}
    \filldraw[violet!100] (1.5,3.5) rectangle (2.5,4.5);
    \filldraw[violet!100] (3.5,3.5) rectangle (4.5,4.5);

    \filldraw[violet!67] (0.5,1.5) rectangle (1.5,2.5);
    \filldraw[violet!67] (2.5,1.5) rectangle (3.5,2.5);
    \filldraw[violet!67] (4.5,1.5) rectangle (5.5,2.5);

    \filldraw[violet!33] (-.5,-.5) rectangle (.5,.5);
    \filldraw[violet!33] (1.5,-.5) rectangle (2.5,.5);
    \filldraw[violet!33] (3.5,-.5) rectangle (4.5,.5);
    \filldraw[violet!33] (5.5,-.5) rectangle (6.5,.5);

    \draw[blue, dashed, very thick] (2,4.55) rectangle (3,5.45);

    \draw[blue, dotted, very thick] (1,2.55) rectangle (2,3.45);
    \draw[blue, dotted, very thick] (3,2.55) rectangle (4,3.45);

    \draw[blue,  very thick] (0,0.55) rectangle (1,1.45);
    \draw[blue,  very thick] (2,0.55) rectangle (3,1.45);
    \draw[blue,  very thick] (4,0.55) rectangle (5,1.45);

    \node[] at (0,0) {$x^6$} ;
    \node[] at (1,0) {$x^5 y$} ;
    \node[] at (2,0) {$x^4 y^2$} ;
    \node[] at (3,0) {$x^3 y^3$} ;
    \node[] at (4,0) {$x^2 y^4$} ;
    \node[] at (5,0) {$x y^5$} ;
    \node[] at (6,0) {$y^6$};

    \node[] at (0.5,1) {$x^5$} ;
    \node[] at (1.5,1) {$x^4 y$} ;
    \node[] at (2.5,1) {$x^3 y^2$} ;
    \node[] at (3.5,1) {$x^2 y^3$} ;
    \node[] at (4.5,1) {$x y^4$} ;
    \node[] at (5.5,1) {$y^5$} ;

    \node[] at (1,2) {$x^4$} ;
    \node[] at (2,2) {$x^3 y$} ;
    \node[] at (3,2) {$x^2 y^2$} ;
    \node[] at (4,2) {$x y^3$} ;
    \node[] at (5,2) {$y^4$} ;

    \node[] at (1.5,3) {$x^3$} ;
    \node[] at (2.5,3) {$x^2 y$} ;
    \node[] at (3.5,3) {$x y^2$} ;
    \node[] at (4.5,3) {$y^3$} ;

    \node[] at (2,4) {$x^2$} ;
    \node[] at (3,4) {$x y$} ;
    \node[] at (4,4) {$y^2$} ;

    \node[] at (2.5,5) {$x$} ;
    \node[] at (3.5,5) {$y$} ;

    \node[] at (3,6) {$1$};

    \end{tikzpicture}
}

    \caption{Figure (a) shows increasing footprints for the moment matrix for $P=2,4,6$ schemes with increasingly lighter shades. The footprints ``support" the corresponding (same shading) even monomial terms in Pascal's triangle, (b). As explained in section \ref{sec:regCellSten}, odd moments of order $P$ are not needed to achieve order $P$ truncation error. In addition, we require at most order $P-1$ monomial terms for the coefficient $\beta$. The cell face $A_{\*i - \*e_x}$ is marked with a red x, and the the dashed, dotted, and solid blue outlines (for $P=2,4,6$, respectively), show the footprint for the modified moment matrix $\*M_{\beta,A_{\*i - \*e_x}}$ in \eqref{eq:modmommat}. The corresponding dashed, dotted, and solid blue boxed monomials in (b) show the order $P-1$ terms that are supported by these footprints. }
    \label{fig:regSten}
\end{figure}

 The vast majority of cells will be regular and will all have the same bilinear stencil, meaning we only have to solve for this stencil once. Since regular cells are squares, the integral of any monomial error term with odd degree over a regular cell is $0$. This means that for the the flux divergence term, in regular cells we can achieve a truncation error of order $P$ using an order $P$ polynomial. We can achieve an order $P$ cell averaged linear term with an order $P-2$ polynomial.
 
 Furthermore, we do not need to calculate all of the Taylor polynomial coefficients up to a given order: if $P$ is even, then the highest order Taylor series coefficients that we need are those $c^{\*q}$ such that $|\*q| = P$ and $q_x, q_y$ are both even. This follows from the fact that our operator does not involve any mixed derivatives: if $|\*q| = P$ for $P$ even, then if $q_x$ is odd, $q_y$ must be as well, so taking derivatives in only one of the dimensions will leave at least one of these powers odd. Although we will get even cross terms from the flux divergence, the even cross terms that arise from combining derivatives of odd order $P$ moments with lower order moments are on the order $P$ and are therefore not needed. These cancellations are typical for centered finite differences, but here they arrive through symmetries in moments and polynomial coefficients. As shown in Figure \ref{fig:regSten}, these coefficients can be supported with a stencil footprint consisting of cells whose centroids are a Manhattan distance of $\frac{P}{2}h$ from the center cell. The columns of the moment matrix $\*M_{u, \*i}$ correspond to all monomials with either $|\*q| < P$ or $|\*q| = P$ and $q_x, q_y$ are both even. Each row is simply the cell-averaged moments $\*m_{\*j}^T$ for each cell $\*j$ in the stencil, and the resulting matrix is square.

 For the matrices $\*M_{\alpha, \*i}$ and $\*M_{\beta, \*i}$ we use the same footprint as $\*M_{u,\*i}$. Construction of $\*s_{\alpha,\*i}$ is then straightforward. The flux divergence term is slightly more complicated. For each of the four faces of the cell $A_{\*i \pm \*e_d}$, if we can use the same flux stencil on each cell's face, then the stencil will guarantee conservation, that is that its contribution to one cell will be the negative of its contribution to its neighbor sharing the same face. For example, consider the face $A_{\*i - \*e_x}$, the left hand vertical face of the cell. For the unit normal, we have $\hat{n}_x = -1$ and $\hat{n}_y =0$, so \eqref{eq:fluxTaylor} reduces to the simpler form:
 \begin{multline}
     \int_{A_{\*i - \*e_x}} \beta \nabla u \cdot \hat{\*n} \ dA  =  \sum_{|\*q|  \leq P}  \lrp{\sum_{|\*r|  \leq (P-|\*q|)}  c_{\beta}^{\*r} \lrp{ - \int_{-\frac{h}{2}}^{\frac{h}{2}} q_x  \lrb{-\frac{h}{2}, y}^{\*q + \*r -\*e_x} dy } } c_{u}^{\*q} + O(h^{P+1}) \ ,
 \end{multline}
meaning the integral in parentheses is only non-zero if $q_x>0$ and $q_y+r_y$ is even. In addition, for accuracy requirements we only need moments with $\abs{\*r + \*q -1} < P$. This means we do not need to calculate any coefficients $\*c_{\beta}^{\*r}$ such that $|\*r|=P$; in other words, an order $P-1$ Taylor approximation to $\beta$ will suffice. Formally, we can multiply the matrix $\*M_{\beta,\*i}$ on the left and right by matrices that zero out the proper columns and rows, giving us the modified moment matrix:
\begin{align}
\label{eq:modmommat}
    \*M_{\beta,A_{\*i - \*e_x}} = \*P_L \*M_{\beta,\*i} \*P_R \ ,
\end{align}
where $\*P_R$ eliminates all columns corresponding to moments with order greater than $P-1$, as well as columns corresponding to moments of order $P-1$ such that $q_y$ is odd. $\*P_L$ eliminates rows that are not necessary to support these moments. We can likewise adjust $\*d_{\beta,\*i}$ and $\*G_{\beta,  \*i}$ to account for these modifications. The contributions from $c_{u}^{\*q}$ are already symmetric about the face because they all represent centered differences. This process results in a flux stencil footprint consisting of all cells whose centroids are a Manhattan distance of $\frac{P-1}{2}h$ from the centroid of the face $A_{\*i - \*e_x}$ (see Figure \ref{fig:regSten}).
However, these simplifications to achieve the minimal stencil footprint rely on symmetry arguments, so they do not apply generally to irregular and cut cells. 

\begin{figure}
\centering

\hspace{0mm}
\subfloat[]{
\begin{tikzpicture}
\filldraw[color=black,fill=yellow!50] (0,0) rectangle (5,5);
    \filldraw[color=red,fill=green!40,very thick] (5,1.5) arc[start angle=270, end angle=180, radius=3.5];
    \fill[color=green!40] (5,1.5) -- (5,5) -- (1.5,5) -- cycle;
    \draw[step=1.0,black,thin] (-1,-1) grid (6,6);

    \node[] at (2.32,2.3) {\footnotesize $V_{+,\*i}$};
    \node[] at (2.7,2.72) {\footnotesize $V_{-,\*i}$};
    \node[] at (1.5,1.5) {\large $\mathcal{N}_{+,\*i}$};
    \node[] at (3.5,3.5) {\large $\mathcal{N}_{-,\*i}$};
    \draw[color=red,very thick] (1.5,4.9) -- (1.5,6);
    \draw[color=red,very thick] (4.9,1.5) -- (6,1.5);

    \draw[color=blue,very thick,dashed] (0,0) rectangle (5,5);
\end{tikzpicture}
}
\subfloat[]{
\begin{tikzpicture}
\filldraw[color=black,fill=yellow!50] (-1,0) rectangle (4,5);

    \filldraw[color=red,fill=green!40,very thick] (5,1.5) arc[start angle=270, end angle=180, radius=3.5];
    \fill[color=green!40] (5,1.5) -- (5,5) -- (1.5,5) -- cycle;

    \filldraw[fill=green!40] (1.5,5) rectangle (6,6);
    \filldraw[fill=green!40] (5,1.5) rectangle (6,6);

    \draw[step=1.0,black,thin] (-1,-1) grid (6,6);

    \node[] at (1.5,2.5) {\footnotesize $V_{+,\*j}$};
    \node[] at (3.5,3.5) {\footnotesize $V_{-,\*k}$};
    \node[] at (.4,1.5) {\large $\mathcal{N}_{+,\*j}$};
    \node[] at (4.5,4.5) {\large $\mathcal{N}_{-,\*k}$};
    
    \draw[color=blue,very thick,dashed] (4,5) -- (4,0) --(-1,0) -- (-1,5) -- (4,5);

    \draw[color=blue,very thick,dashed] (1,1) rectangle (6,6);
    
    \draw[color=red,very thick] (5,1.5) arc[start angle=270, end angle=180, radius=3.5];
    
    \draw[color=red,very thick] (1.5,4.9) -- (1.5,6.022);
    \draw[color=red,very thick] (4.9,1.5) -- (6.022,1.5);
\end{tikzpicture}
}
    \caption{Figure (a) shows neighborhoods $\mathcal{N}_{\pm,\*i}$ used to construct interpolation matrices for order $P=2$, around cut cell $\*i$, which contains two volumes $V_{\pm,\*i}$ bordering the interface. Figure (b) shows neighborhoods surrounding full cells $\*j$ and $\*k$ that don't contain the interface, but are ``irregular,'' meaning the regular stencil for order $P=2$ would be inconsistent ($O$th-order).}
    \label{fig:irregSten}
\end{figure}
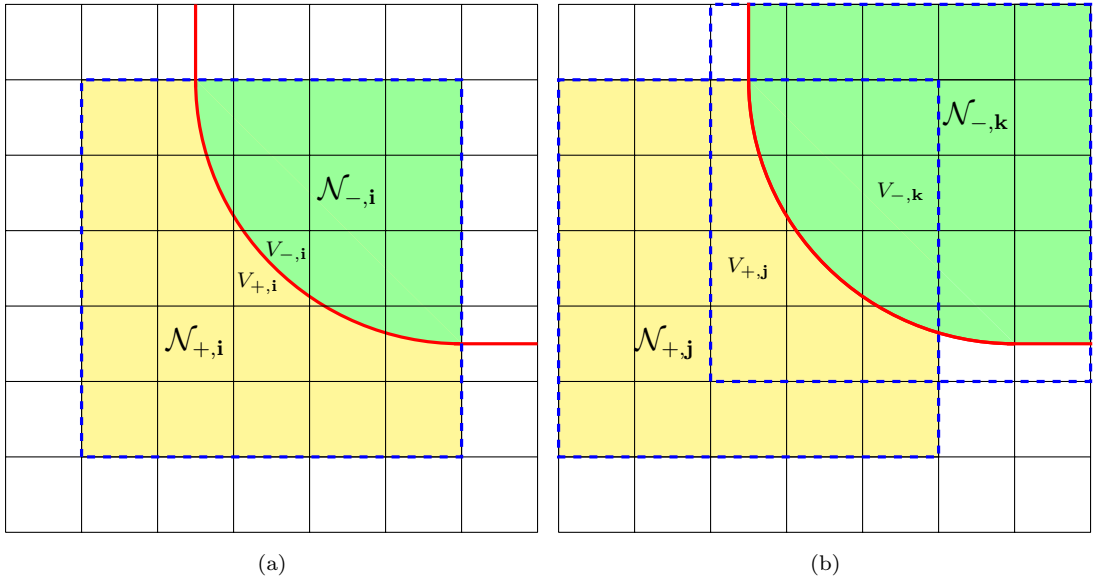

\subsubsection{Irregular Cells}
Let cell $\*i$ be an irregular cell in phase $p$, meaning its regular cell footprint contains at least one cell which is intersected by the EB. This invalidates the truncation error analysis for the regular cell stencil, so we adopt a more general method for construction of moment matrices for irregular cells. Let $\mathcal{N}_{p, \*i}$ be a neighborhood of cells in phase $p$ around cell $\*i$ (see Figure \ref{fig:irregSten}(b) \,). Our data vector $\*d_{u,p, \*i}$ will consist of cell averaged values $\IP{u}_{\*j}$ for each cell $\*j \in \mathcal{N}_{p, \*i}$. Each row of the corresponding moment matrix $\*M_{u,p,\*i}$ will simply be $\*m_{\*j}^T$, the row of cell averaged volume moments for each cell $\*j$. For irregular cells, the moment matrices $\*M_{\alpha, p,\*i}$ and $\*M_{\beta, p,\*i}$ are identical.  We form one of each per cell and use a different $\*G_{\beta, p, \*i \pm \*e_d}$ for each face to create a stencil for the flux integral along a face. The columns of these moment matrices consist of all moments up to order $P$.

\subsubsection{Cut Cells}
In cut cells, the moment matrix is additionally used to enforce jump conditions. When interpreted as an overdetermined system, our interpolating polynomials are being constrained to satisfy interface matching conditions. Alternatively as an underdetermined system, we use jump conditions as data to cancel truncation error terms. It is necessary to enforce these jump conditions so that our discrete operator is not degenerate. 

Let $\mathcal{N}_{p,\*i}$ be two neighborhoods of cells in their respective phases around the cut cell $\*i$. See Figure \ref{fig:irregSten} (a). As with irregular cells, we can form two moment matrices $\*M_{u,p,\*i}$ and data vectors $\*d_{u,p,\*i}$ of cell averaged values. For each cut cell $\*j$ in $\mathcal{N}_{+,\*i} \cup \mathcal{N}_{-,\*i}$, we want to enforce the two jump conditions \eqref{eq:uJumpD} and \eqref{eq:duJumpD}. Expressing these in terms of moments and Taylor coefficients, we have 
\begin{align}
    \int_{A_{B,\*j}} u^+ - u^- \ dA &= \sum_{|\*q| \leq P} (c_{u,+}^{\*q} - c_{u,-}^{\*q}) m_{B, \*j}^{\*q} + O(h^{P+2})  \ ,
\end{align}
and for the jump in the flux:
\begin{multline}
    \int_{A_{B,\*j}} (\beta^+ \nabla u^+ - \beta^- \nabla u^-) \cdot \hat{\*n} \ d A = 
    \sum_{|\*r|,|\*q|  \leq P}  \lrp{c_{\beta,+}^{\*r} c_{u,+}^{\*q} - c_{\beta,-}^{\*r} c_{u,-}^{\*q} }  \lrb{ q_x m_{B, \*j, x}^{\*q + \*r -\*e_x} + q_y m_{B, \*j, y}^{\*q + \*r -\*e_y} } + O(h^{P+1}) \ .
\end{multline}
These expressions are both linear in the coefficients $c_{u,p}^{\*q}$, and they couple both sets of coefficients by interpolating $u^+$ and $u^-$ simultaneously. The resulting matrices $\*M_{J,+, \*i}$ and $\*M_{J,-, \*i}$ have rows with the jump condition expressions and columns corresponding to $\*c_{u,+}$ and $\*c_{u,-}$, respectively. The vector $\*d_{J,\*i}$ consists of the given jump condition data. Finally we form the moment matrix $\*M_{\*i}$ which is defined as:
\begin{align}
    \*M_{\*i} = \begin{bmatrix} \*M_{u,+,\*i} & 0\\
    0 & \*M_{u,-,\*i} \\
    \*M_{J,+, \*i} & \*M_{J,-, \*i} 
    \end{bmatrix} \ ,
\end{align}
and we solve for both sets of coefficients simultaneously:
\begin{align}
    \begin{bmatrix} \*c_{u,+} \\ 
    \*c_{u,-}  \end{bmatrix} &= \*M_{\*i}^{\dagger} \begin{bmatrix}
    \*d_{u,+,\*i} \\
    \*d_{u,-,\*i} \\
    \*d_{J,\*i}
    \end{bmatrix} 
    = \begin{bmatrix}
    (\*M_{\*i}^{\dagger})_V & (\*M_{\*i}^{\dagger} )_J
    \end{bmatrix} 
    \begin{bmatrix}
    \*d_{u,+,\*i} \\
    \*d_{u,-,\*i} \\
    \*d_{J,\*i}
    \end{bmatrix} = (\*M_{\*i}^{\dagger})_V \begin{bmatrix}
    \*d_{u,+,\*i} \\
    \*d_{u,-,\*i} 
    \end{bmatrix} + (\*M_{\*i}^{\dagger} )_J \*d_{J,\*i} \ .
\end{align}
If values in $\*d_{J,\*i}$ are nonzero, then forming a stencil using these coefficients will result in adding a scalar to the right hand side $\IP{f}_{p,\*i}$. For example, for the linear term we would have:
\begin{align}
    \*s_{\alpha, p, \*i}^T &= \lrb{\*d_{\alpha, p, \*i}^T (\*M_{\alpha,p, \*i}^{\dagger})^T \*G_{\alpha, p, \*i}^T (\*M_{\*i}^{\dagger})_V} + \lrb{\*d_{\alpha, p, \*i}^T (\*M_{\alpha,p, \*i}^{\dagger})^T \*G_{\alpha, p, \*i}^T (\*M_{\*i}^{\dagger} )_J \*d_{J,\*i}} \ ,
\end{align}
where the second term in brackets is a scalar. To form the moment matrices $\*M_{\alpha, p,\*i}$ and $\*M_{\beta, p,\*i}$, we can just use the same neighborhoods as the $\*M_{u,p,\*i}$; we do not need to couple these systems because there are no external constraints on the jumps in coefficients. 

\subsubsection{Conservation}
For each cut cell and irregular cell, we have shown how to obtain a stencil for $\int_A \beta \nabla u \cdot \hat{\*n} \ dA$ on each face $A$ of any cell. However, in order to have a conservative method, we must have only one flux stencil per non-EB face. A simple solution to this problem is to average the flux stencils between pairs of neighboring irregular and cut cells. Although this results in larger stencils, it has the advantage of coupling a layer of irregular cells to the interface jump conditions. This is because the irregular cell stencils that border cut cells will share stencil information with cut cells, which incorporate interface jump conditions. We reiterate that this is not a significant issue because the density of the linear system is dominated by the size of the regular cell stencil. For regular cells we do not have to average with neighbors because our flux stencils were created individually for each face and are symmetric about that face. At a cell face which is shared between an irregular and regular cell, we use the regular cell flux stencil.

\subsection{Neighborhood Selection and Weighting}
In general, neighborhoods need to be chosen so that resulting moment matrices are overdetermined. So as to not perform some sort of search based on local geometry, we opt to make the neighborhood sufficiently large to accommodate a reasonably smooth geometry. For any irregular or cut cell in phase $p$, we let $\mathcal{N}_{p, \*i}$ be those cells in phase $p$ that lie in the the square of cells with side length $2P +1$ surrounding cell $i$. See Figure \ref{fig:irregSten}. As is done in \cite{Dharshi} and \cite{nate}, we employ a weighted least-squares approach to force stencil weights to decay with distance faster than the growth of the highest polynomial term. To each piece of data in $\*d$ we assign a weight that is inversely related to its distance from the centroid of volume $\Vpi$. If $\delta_{\*j}$ is this distance, then the corresponding row of the moment matrix and $\*d$ are multiplied by $w_{\*j}$, where
\begin{align}
    w_{\*j} = \frac{1}{(1+\delta_{\*j})^{P+1}} \ ,
\end{align}
where $P$ is the order of the scheme.
This forms a diagonal weight matrix $\*W$, which gives us the weighted least-squares solution:
\begin{align}
    \*c_{u} = (\*W \*M)^{\dagger} \*W \*d \ .
\end{align}
This weighting does not affect the truncation error: with some matrix algebra we can see that this weighted least-squares solution is equivalent to a change of basis in the undetermined formulation:
\begin{align}
    \*s &= \*W \lrp{\*M^T \*W}^{\dagger} \*g \implies \\
    \*M^T \*W \*W^{-1} \*s &= \*g \ .
\end{align}
Since $\*M^T \*s = \*g$ still holds, the truncation error is unaffected by the weighting. This has proved to be an effective tool for controlling the spectrum and conditioning of the discrete operator \cite{Dharshi}. We compute the pseudoinverse using the SVD algorithm in 
\emph{LAPACK} \cite{LAPACK}.

\subsection{Solver and Software Implementation}
We assemble the stencils to form the linear system
\begin{align}
    \*L \*u = \*f + \*r \ ,
\end{align}
where $\*r$ represents the contribution to the right hand side from the jump conditions. Following \cite{Dharshi}, we precondition this system by left multiplying with the diagonal matrix whose $\lrp{i,i}$ entry is $\frac{|\Vi|}{h^2}$; i.e. we multiply each row by the volume fraction of that cell. This simple preconditioner eliminates the volume scaling associated with very small volume fractions.

We solve the linear system using Krylov subspace methods and preconditioners provided by the \emph{PETSc} library \cite{petsc-user-ref}, \cite{petsc-efficient}. Since the linear system is non-symmetric, we use BiCG-Stab or GMRES. We have experimented with the \emph{PETSc} algebraic multigrid and block Jacobi preconditioners. One of these options is typically sufficient, but if they fail, we use the direct solver \emph{SuperLU} \cite{lidemmel03}. In future research we will develop a geometric multigrid preconditioner similar to that in \cite{Crockett} or \cite{Devendran2015AHM}. Our method is well-suited for geometric multigrid because the Taylor series formulation makes interpolation straightforward. However, we emphasize that an efficient solver is not the focus of this particular paper. The algorithm is implemented using the \emph{Chombo} software library \cite{Chombo}, which allows for large-scale parallelization of the algorithm. Visualizations are created using \emph{VisIt} \cite{visit}.



\section{Numerical Tests}
We validate our method with a series of numerical tests. The goals of this section are to:
\begin{enumerate}
    \item Validate the truncation error analysis and measure the solution error,
    \item Demonstrate consistent accuracy for problems with large coefficient and solution jumps, and
    \item Demonstrate convergence on non-trivial interface geometries.
\end{enumerate}
Although our scheme is designed to be have arbitrary order of accuracy, we have evaluated it for just $P \in \set{2,4,6}$.  We measure the error $\IPi{e}$ as discrete cell averages, and evaluate it using discrete $L^p$ norms:
\begin{align}
    \Norm{e}_1 &= \int_{\Omega} |e| dV = \sum_{p,\*i} \abs{\IP{e}_{p, \*i}} |V_{p,\*i}| \\
    \Norm{e}_{\infty} &= \max_{p,\*i}  \IP{e}_{p, \*i} \ .
\end{align}

\begin{figure}
    \centering
    \subfloat[]{\includegraphics[width=120mm]{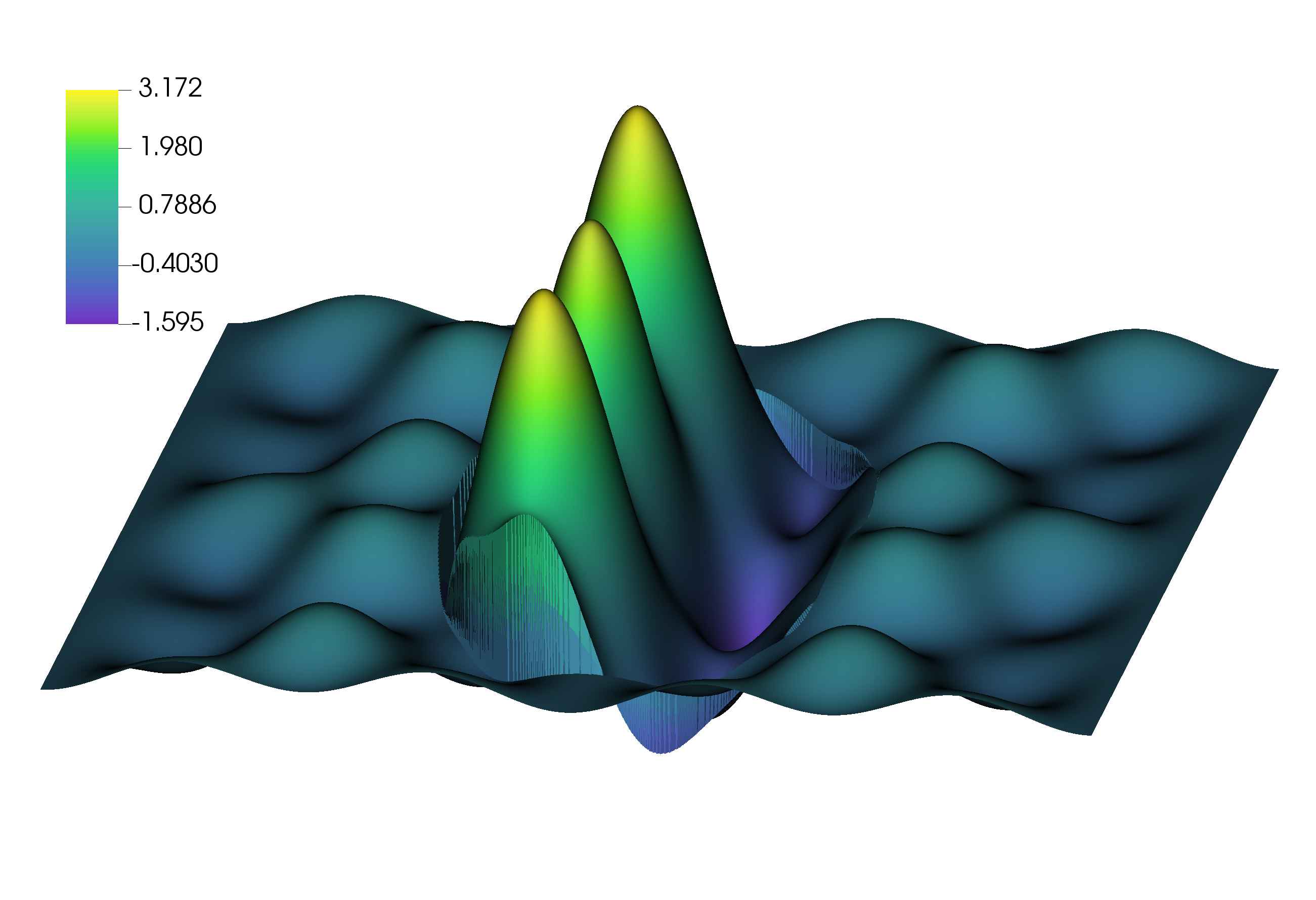}}
    \hspace{0mm}
    \subfloat[]{\includegraphics[width=85mm]{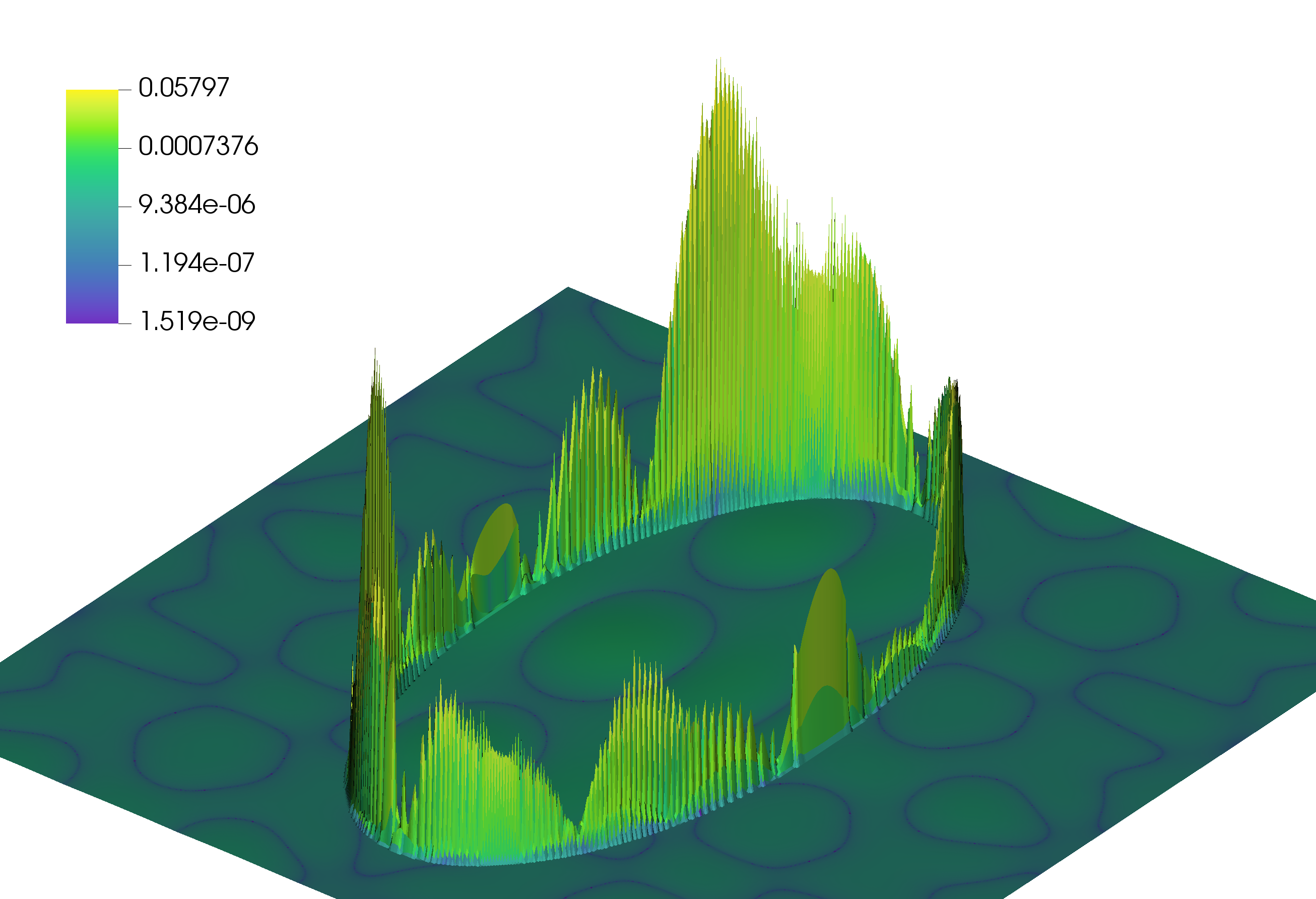}}
    \subfloat[]{\includegraphics[width=85mm]{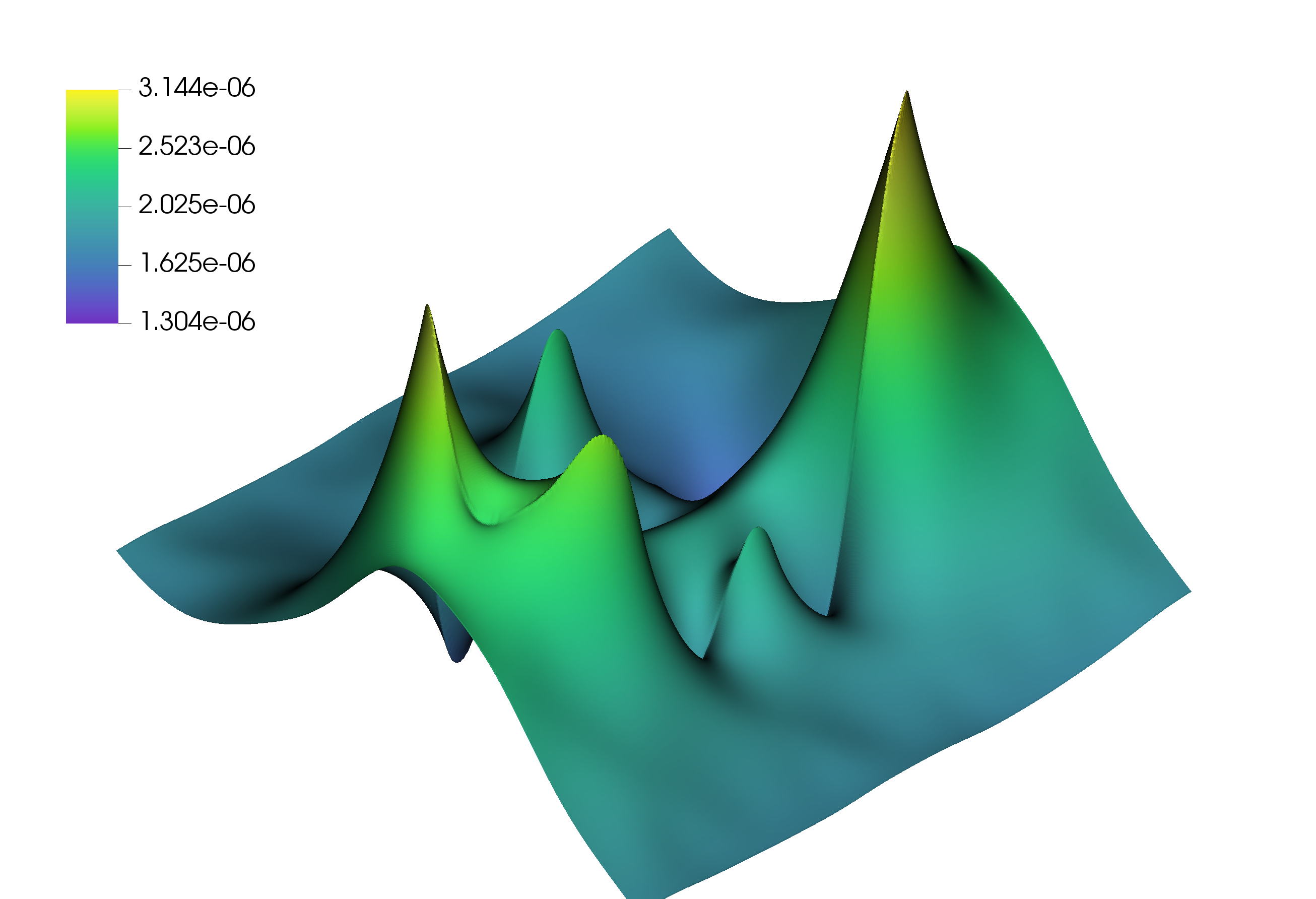}}
    \caption{Results of the ellipse boundary tests for $P=4$, $n=512$. Plots of the (a) exact solution, (b) absolute value of truncation error (log scale), and (c) absolute value of solution error (log scale). Note that the truncation error is concentrated at the interface, and its influence on the solution error is much smoother after inverting the elliptic operator.}
    \label{fig:Error}
\end{figure}

\begin{figure}
\centering
\subfloat[]{
\begin{tikzpicture}[scale=.9]
\begin{axis}[
    xmode = log,
    ymode = log,
    title={Truncation Error, $L^{\infty}$ Norm},
    xlabel={$n$},
    ylabel={Error},
    xmin=32, xmax=512,
    ymin=1e-6, ymax=1e3,
    xtick={32,64,128,256,512},
    xticklabels={32,64,128,256,512},
    ytick={1e-6,1e-4,1e-2,1e0,1e2},
    legend pos=south west,
    ymajorgrids=true,
    grid style=dashed,
]

\addplot[
    color=blue,
    mark=square,
    ]
    coordinates {
    (32,1.264e+02)(64,5.157e+00)(128,2.138e-01) (256,6.695e-03)(512,2.023e-04)
    };
    \addlegendentry{$P=6$ }
\addplot[
    color=blue,
    dashed,
    ]
    coordinates {
    (32,1.264e+02)(512,1.2054e-04)
    };
    \addlegendentry{$O(h^5)$}
\addplot[
    color=red,
    mark=square,
    ]
    coordinates {
    (32,1.443e+02)(64,2.391e+01)(128,3.851e+00) (256,4.804e-01)(512,5.797e-02)
    };
    \addlegendentry{$P=4$ }
\addplot[
    color=red,
    dashed,
    ]
    coordinates {
    (32,1.443e+02)(512,0.1417)
    };
    \addlegendentry{$O(h^3)$}
\addplot[
    color=cyan,
    mark=square,
    ]
    coordinates {
    (32,1.837e+02)(64,8.299e+01)(128,4.240e+01) (256,2.229e+01)(512,1.055e+01)
    };
    \addlegendentry{$P=2$ }
\addplot[
    color=cyan,
    dashed,
    ]
    coordinates {
    (32,1.837e+02)(512,11.4812)
    };
    \addlegendentry{$O(h)$}

\end{axis}
\end{tikzpicture}
}
\subfloat[]{
\begin{tikzpicture}[scale=.9]
\begin{axis}[
    xmode = log,
    ymode = log,
    title={Truncation Error, $L^{1}$ Norm},
    xlabel={$n$},
    ylabel={Error},
    xmin=32, xmax=512,
    ymin=1e-6, ymax=1e3,
    xtick={32,64,128,256,512},
    xticklabels={32,64,128,256,512},
    ytick={1e-6,1e-4,1e-2,1e0,1e2},
    legend pos=south west,
    ymajorgrids=true,
    grid style=dashed,
]

\addplot[
    color=blue,
    mark=square,
    ]
    coordinates {
    (32,1.829e+01)(64,3.385e-01 )(128,5.344e-03) (256,7.974e-05)(512,1.231e-06)
    };
    \addlegendentry{$P=6$ }
\addplot[
    color=blue,
    dashed,
    ]
    coordinates {
    (32,1.829e+01)(512,1.0902e-06)
    };
    \addlegendentry{$O(h^6)$}
\addplot[
    color=red,
    mark=square,
    ]
    coordinates {
    (32,2.774e+01)(64,1.676e+00)(128,1.048e-01) (256,6.343e-03)(512,3.927e-04)
    };
    \addlegendentry{$P=4$ }
\addplot[
    color=red,
    dashed,
    ]
    coordinates {
    (32,2.774e+01)(512,4.2328e-04)
    };
    \addlegendentry{$O(h^4)$}
\addplot[
    color=cyan,
    mark=square,
    ]
    coordinates {
    (32,7.646e+01)(64,1.831e+01)(128,4.512e+00) (256,1.106e+00)(512,2.756e-01)
    };
    \addlegendentry{$P=2$ }
\addplot[
    color=cyan,
    dashed,
    ]
    coordinates {
    (32,7.646e+01)(512,0.2987)
    };
    \addlegendentry{$O(h^2)$}

\end{axis}
\end{tikzpicture}
}
\hspace{0mm}
\subfloat[]{
\begin{tikzpicture}[scale=.9]
\begin{axis}[
    xmode = log,
    ymode = log,
    title={Solution Error, $L^{\infty}$ Norm},
    xlabel={$n$},
    ylabel={Error},
    xmin=32, xmax=512,
    ymin=1e-8, ymax=1e0,
    xtick={32,64,128,256,512},
    xticklabels={32,64,128,256,512},
    ytick={1e-8,1e-6,1e-4,1e-2,1e0},
    legend pos=south west,
    ymajorgrids=true,
    grid style=dashed,
]

\addplot[
    color=blue,
    mark=square,
    ]
    coordinates {
    (32,9.505e-02)(64,2.291e-03)(128,4.132e-05) (256,6.802e-07)(512,1.075e-08)
    };
    \addlegendentry{$P=6$ }
\addplot[
    color=blue,
    dashed,
    ]
    coordinates {
    (32,9.505e-02)(512,5.6654e-09)
    };
    \addlegendentry{$O(h^6)$}
\addplot[
    color=red,
    mark=square,
    ]
    coordinates {
    (32,1.360e-01)(64,1.179e-02)(128,7.838e-04) (256,5.025e-05)(512,3.146e-06)
    };
    \addlegendentry{$P=4$ }
\addplot[
    color=red,
    dashed,
    ]
    coordinates {
    (32,1.360e-01)(512,2.0752e-06)
    };
    \addlegendentry{$O(h^4)$}
\addplot[
    color=cyan,
    mark=square,
    ]
    coordinates {
    (32,2.456e-01)(64,5.312e-02)(128,1.214e-02) (256,2.990e-03)(512,7.539e-04)
    };
    \addlegendentry{$P=2$ }
\addplot[
    color=cyan,
    dashed,
    ]
    coordinates {
    (32,2.456e-01)(512,9.5938e-04)
    };
    \addlegendentry{$O(h^2)$}

\end{axis}
\end{tikzpicture}
}
\subfloat[]{
\begin{tikzpicture}[scale=.9]
\begin{axis}[
    xmode = log,
    ymode = log,
    title={Solution Error, $L^{1}$ Norm},
    xlabel={$n$},
    ylabel={Error},
    xmin=32, xmax=512,
    ymin=1e-8, ymax=1e0,
    xtick={32,64,128,256,512},
    xticklabels={32,64,128,256,512},
    ytick={1e-8,1e-6,1e-4,1e-2,1e0},
    legend pos=south west,
    ymajorgrids=true,
    grid style=dashed,
]

\addplot[
    color=blue,
    mark=square,
    ]
    coordinates {
    (32,2.575e-01)(64,5.910e-03)(128,1.048e-04) (256,1.755e-06)(512,2.811e-08 )
    };
    \addlegendentry{$P=6$ }
\addplot[
    color=blue,
    dashed,
    ]
    coordinates {
    (32,2.575e-01)(512,1.5348e-08)
    };
    \addlegendentry{$O(h^6)$}
\addplot[
    color=red,
    mark=square,
    ]
    coordinates {
    (32,3.533e-01)(64,2.962e-02)(128,1.934e-03) (256,1.262e-04)(512,7.987e-06)
    };
    \addlegendentry{$P=4$ }
\addplot[
    color=red,
    dashed,
    ]
    coordinates {
    (32,3.533e-01)(512,5.3909e-06)
    };
    \addlegendentry{$O(h^4)$}
\addplot[
    color=cyan,
    mark=square,
    ]
    coordinates {
    (32,1.613e-01)(64,4.237e-02)(128,1.019e-02) (256,2.526e-03)(512,6.371e-04)
    };
    \addlegendentry{$P=2$ }
\addplot[
    color=cyan,
    dashed,
    ]
    coordinates {
    (32,1.613e-01)(512,6.3008e-04)
    };
    \addlegendentry{$O(h^2)$}

\end{axis}
\end{tikzpicture}
}
\caption{Truncation error $L^\infty$ norm (a), $L^1$ norm (b), and solution error norms (c) and (d), respectively, for the ellipse geometry. We observe the expected rate of convergence for both quantities in both norms. The truncation error is dominated by behavior at the interface, which because it is codomension-1 smaller, achieves one order higher in $L^1$ norm (see Figure \ref{fig:Error}). }
\label{fig:truncError}
\end{figure}
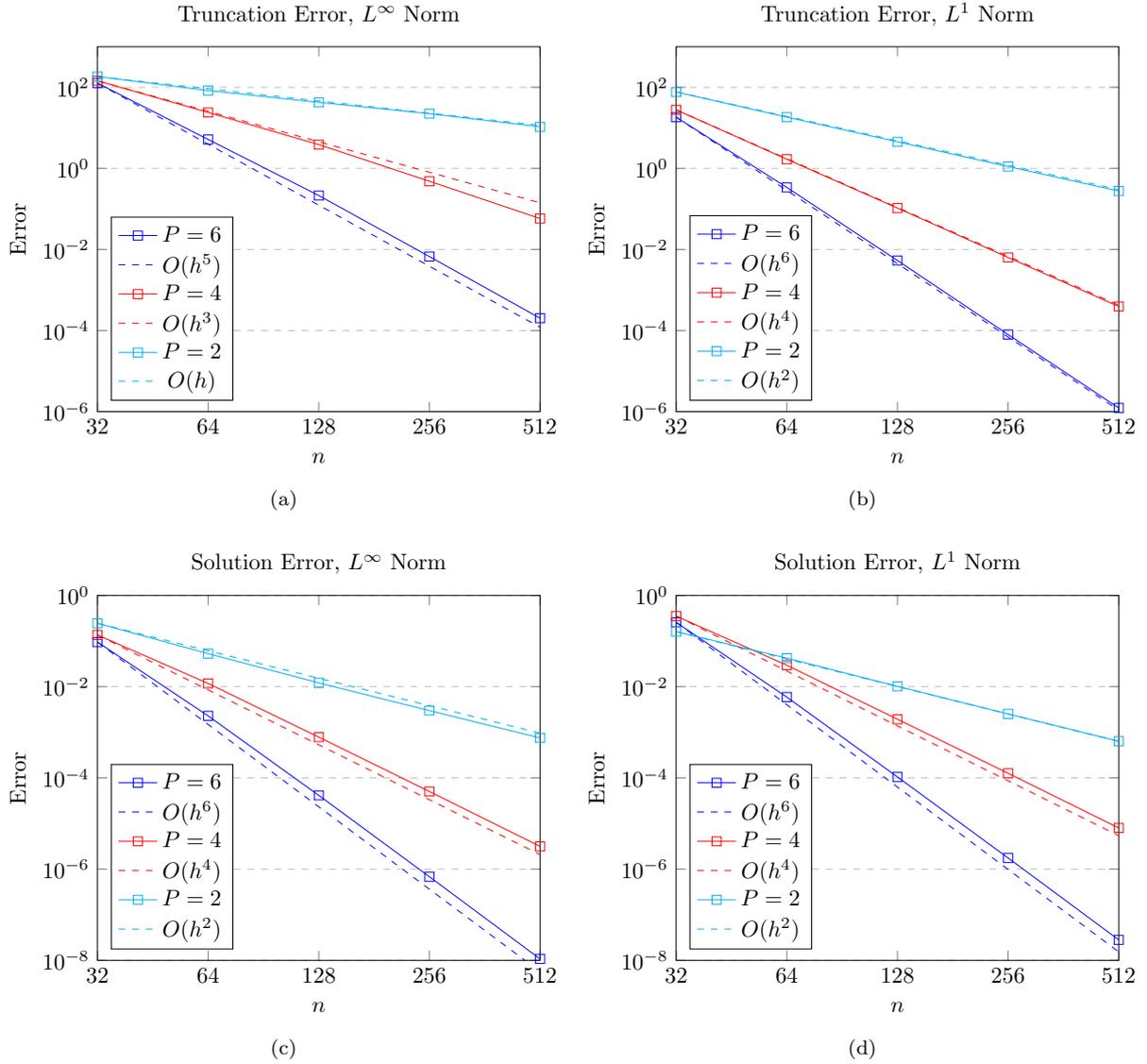

\subsection{Truncation and Solution Error Validation}
\label{sec:ellipse_test}
Our physical domain for all tests is $\Omega = [-1,1]^2$, intersected by some interface $\Gamma$. Our first interface $\Gamma$ is an ellipse with major axis of length $\frac{\pi}{4}$ and minor axis of length $\frac{\pi}{8}$. The superscript $+$ refers to quantities enclosed by the interface and the superscript $-$ refers to quantities on the exterior of the interface. We test our discretization using the method of manufactured solutions, such that $\alpha^{\pm}, \beta^{\pm}, u^{\pm}$ are all constructed as linear combinations of periodic functions $p_{k_x, k_y}$ on the square:
\begin{align}
    p_{k_x, k_y} = \cos^2 \lrp{\pi k_x x } \sin^2 \lrp{\pi k_y y } \ ,
\end{align}
and let $k_x, k_y \in \set{-2,-1,1,2}$, creating 20 total basis functions. We randomly generate different sets of coefficients $c_{k_x, k_y} \in [-1,1]$ for each of the six functions $\alpha^{\pm}, \beta^{\pm}, u^{\pm}$. An exact $f^{\pm}$ is formed by applying the exact differential operator to $u^{\pm}$. The variable coefficient fields are offset by a positive constant so that they are nonnegative everywhere. Figure \ref{fig:Error}(a) shows the solution for one particular example.

The global truncation error $\*e$ and solution error $\*t$ are defined as
\begin{align}
    \*e &= \*u^e - \*u \\
    \*t &= \lrp{\*L \*u^e - \*r  }- \*L^e \*u^e \ ,
\end{align}
where $\*u^e$ is the exact solution and $\*L^e$ is the exact operator, so that the solution error satisfies the equation
\begin{align}
    \*L \*e = \*t \label{eq:errorEq} \ .
\end{align}
Note that \eqref{eq:errorEq} implies homogeneous jump conditions on the error; there is no contribution to the right-hand-side of this system from jump conditions.

Our truncation error analysis predicts an order $P-1$ truncation error in cut and irregular cells and an order $P$ truncation error in regular cells. In Figure \ref{fig:Error}(b) we see that truncation error (for $P=4$) is almost entirely concentrated in cut and irregular cells. This is reflected in Figure \ref{fig:truncError}(a), as the max norm of the truncation error converges at order $P-1$, the expected rate for cells at and near the interface. However, the number of cut and irregular cells is of order $h^{-1}$ because it is a codimension one smaller region, while the number of regular cells is scales like $h^{-2}$. Therefore for the $L^1$ norm we have:
\begin{align}
    \Norm{t}_1 &= \sum_{\*i \in \Omega_C \cup \Omega_I} \abs{\IP{t}_{p,\*i}} |V_{p,\*i}| + \sum_{\*i \in \Omega_R} \abs{\IP{t}_{p,\*i}} |V_{p,\*i}|  = \sum_{\*i \in \Omega_C \cup \Omega_I} O(h^{P+1}) + \sum_{\*i \in \Omega_R}  O(h^{P+2}) \\
    &= O(h^{-1}) O(h^{P+1}) + O(h^{-2}) O(h^{P+2}) = O(h^P) \ .
\end{align}
This is confirmed in \ref{fig:truncError} (b); we see clean order $P$ convergence for the $L^1$ norm of the truncation error. Given that the truncation error at the interface is orders of magnitude greater than truncation error elsewhere, the $L^1$ norm of the truncation error is also dominated by behavior at the interface.

Although we do not have an analytical bound on $\Norm{\*L^{-1}}$, based on the analysis in \cite{Hans_thesis} and the results in \cite{Dharshi}, \cite{Crockett}, \cite{gibou} and others, we expect the solution error to converge at order $P$ in both norms. This behavior is shown in Figures \ref{fig:Error}(b) and (c); the solution error is roughly of the same order of magnitude everywhere in the domain. We plan on further analysis to explore the combined effects of the homogeneous jump conditions imposed on the error equation and the regularity of the elliptic operator, but the empirical results demonstrate the desired convergence rates.

\FloatBarrier

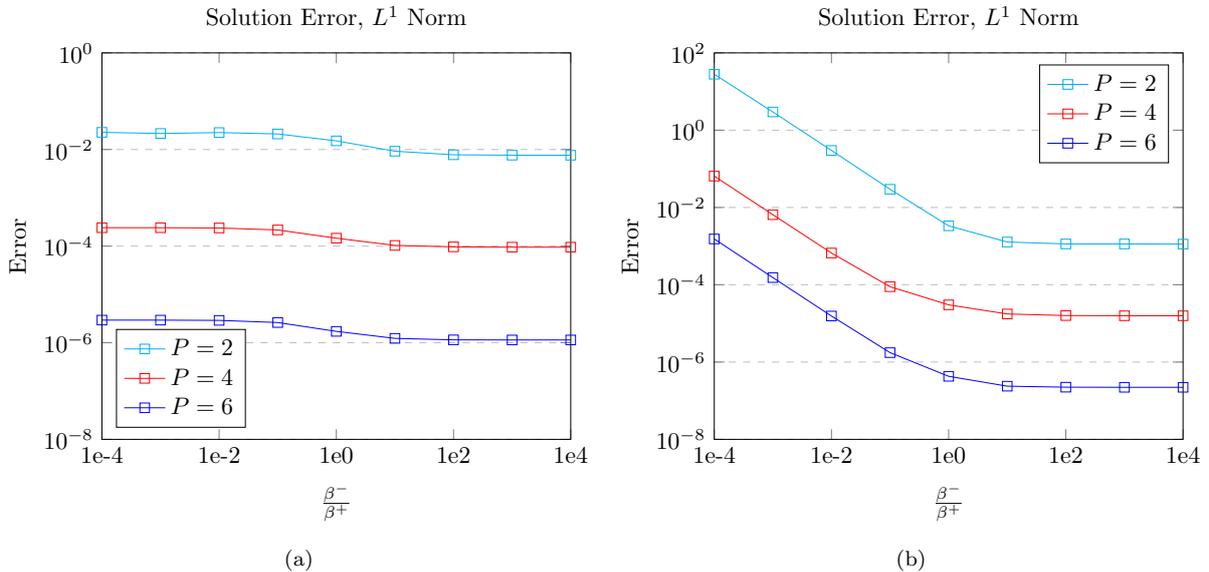
\begin{figure}
\centering

\subfloat[]{
\begin{tikzpicture}[scale=.9]
    \begin{axis}[
    ymode = log,
    xmode = log,
    title={Solution Error, $L^1$ Norm},
    xlabel={$\frac{\beta^-}{\beta^+}$},
    ylabel={Error},
    xmin=1e-4, xmax=1e4,
    ymin=1e-8, ymax=1e-0,
    xtick={1e-4, 1e-2,  1e0,  1e2,  1e4},
    xticklabels={1e-4, 1e-2,  1e0,  1e2,  1e4},
    ytick={1e-8, 1e-6, 1e-4, 1e-2,1e0},
    legend pos=south west,
    ymajorgrids=true,
    grid style=dashed,
]

\addplot[
    color=cyan,
    mark=square,
    ]
    coordinates {
    (1e-4, 2.266e-02) (1e-3,2.145e-02 ) (1e-2,2.242e-02 )(1e-1, 2.098e-02) (1e0,1.510e-02 )(1e1,9.217e-03 )(1e2,7.781e-03 )(1e3,7.630e-03 )(1e4,7.613e-03 )
    };
    \addlegendentry{$P=2$}

\addplot[
    color=red,
    mark=square,
    ]
    coordinates {
    (1e-4, 2.405e-04) (1e-3,2.402e-04 ) (1e-2,2.374e-04 )(1e-1,2.160e-04) (1e0,1.459e-04 )(1e1,1.033e-04 )(1e2,9.649e-05 )(1e3,9.508e-05 )(1e4,9.502e-05 )
    };
    \addlegendentry{$P=4$}
\addplot[
    color=blue,
    mark=square,
    ]
    coordinates {
    (1e-4, 2.941e-06) (1e-3,2.937e-06 ) (1e-2,2.902e-06 )(1e-1,2.620e-06) (1e0,1.719e-06 )(1e1,1.231e-06 )(1e2,1.147e-06 )(1e3, 1.141e-06 )(1e4,1.141e-06 )
    };
    \addlegendentry{$P=6$}

\end{axis}
\end{tikzpicture}
}
\subfloat[]{
\begin{tikzpicture}[scale=.9]
    \begin{axis}[
    ymode = log,
    xmode = log,
    title={Solution Error, $L^1$ Norm},
    xlabel={$\frac{\beta^-}{\beta^+}$},
    ylabel={Error},
    xmin=1e-4, xmax=1e4,
    ymin=1e-8, ymax=1e2,
    xtick={1e-4, 1e-2,  1e0,  1e2,  1e4},
    xticklabels={1e-4, 1e-2,  1e0,  1e2,  1e4},
    ytick={1e-8, 1e-6, 1e-4, 1e-2,1e0,1e2},
    legend pos=north east,
    ymajorgrids=true,
    grid style=dashed,
]

\addplot[
    color=cyan,
    mark=square,
    ]
    coordinates {
    (1e-4, 2.800e+01) (1e-3,2.962e+00 ) (1e-2,2.989e-01 )(1e-1, 2.957e-02) (1e0,3.348e-03 )(1e1,1.276e-03 )(1e2,1.139e-03 )(1e3,1.137e-03 )(1e4,1.129e-03 )
    };
    \addlegendentry{$P=2$}

\addplot[
    color=red,
    mark=square,
    ]
    coordinates {
    (1e-4, 6.481e-02) (1e-3,6.496e-03 ) (1e-2,6.625e-04 )(1e-1,8.904e-05) (1e0,3.018e-05 )(1e1,1.754e-05 )(1e2,1.596e-05 )(1e3,1.577e-05 )(1e4,1.576e-05 )
    };
    \addlegendentry{$P=4$}
\addplot[
    color=blue,
    mark=square,
    ]
    coordinates {
    (1e-4, 1.531e-03 ) (1e-3,1.533e-04 ) (1e-2,1.552e-05 )(1e-1,1.752e-06) (1e0,4.276e-07 )(1e1,2.375e-07 )(1e2,2.237e-07  )(1e3, 2.221e-07 )(1e4,2.221e-07 )
    };
    \addlegendentry{$P=6$}

\end{axis}
\end{tikzpicture}
}
\caption{Solution error for discontinuous diffusion coefficient tests. Figure (a) is for the cosine geometry, whereas (b) is the annulus geometry (we only plot $L^1$ error because the $L^{\infty}$ norms behave nearly identically). }
\label{fig:coefTestError}
\end{figure}
\begin{figure}
    \centering
    \subfloat[]{
    \includegraphics[width=70mm]{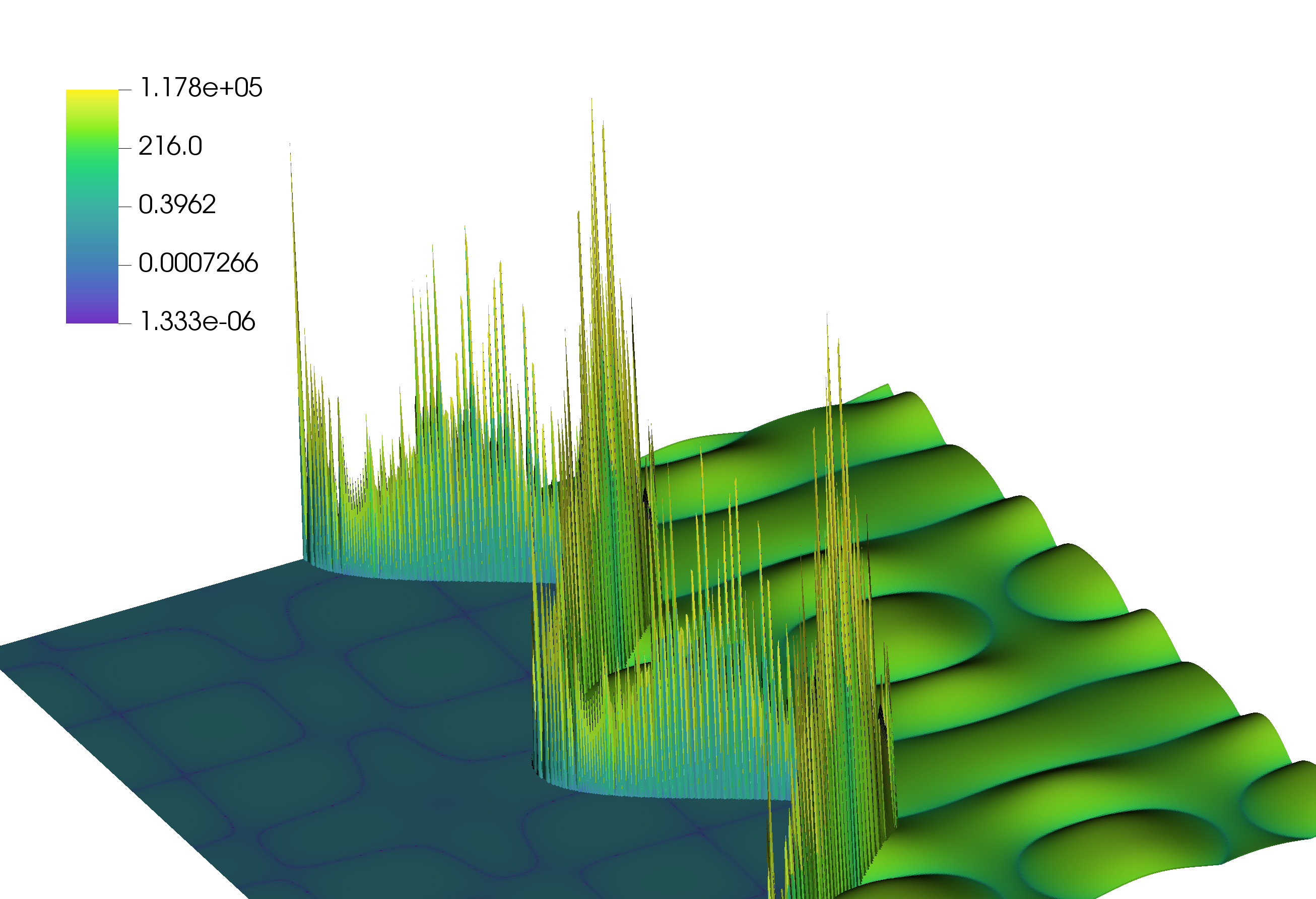}}
    \subfloat[]{
    \includegraphics[width=70mm]{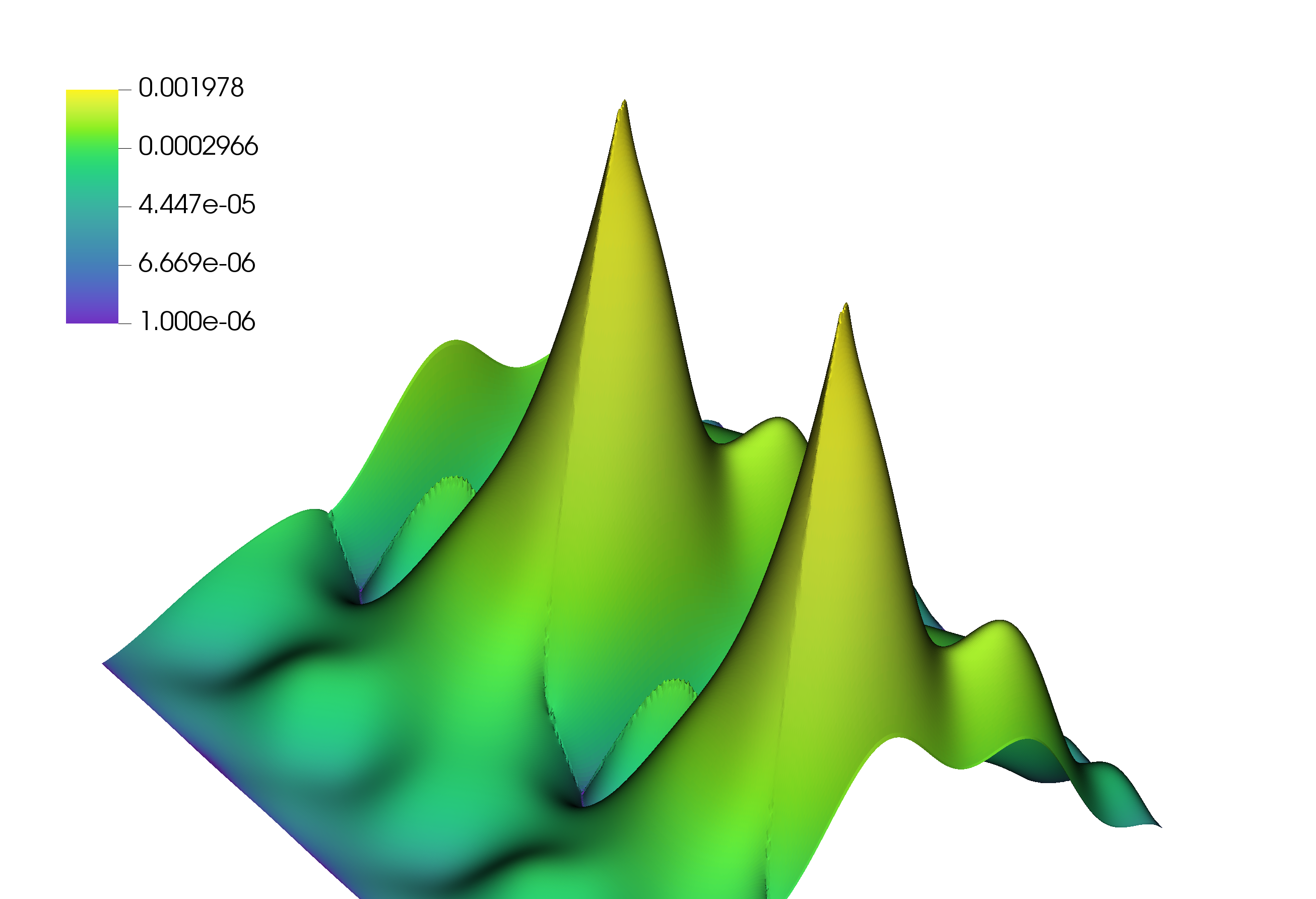}
    }
    \hspace{0 mm}
    \subfloat[]{
    \includegraphics[width=70mm]{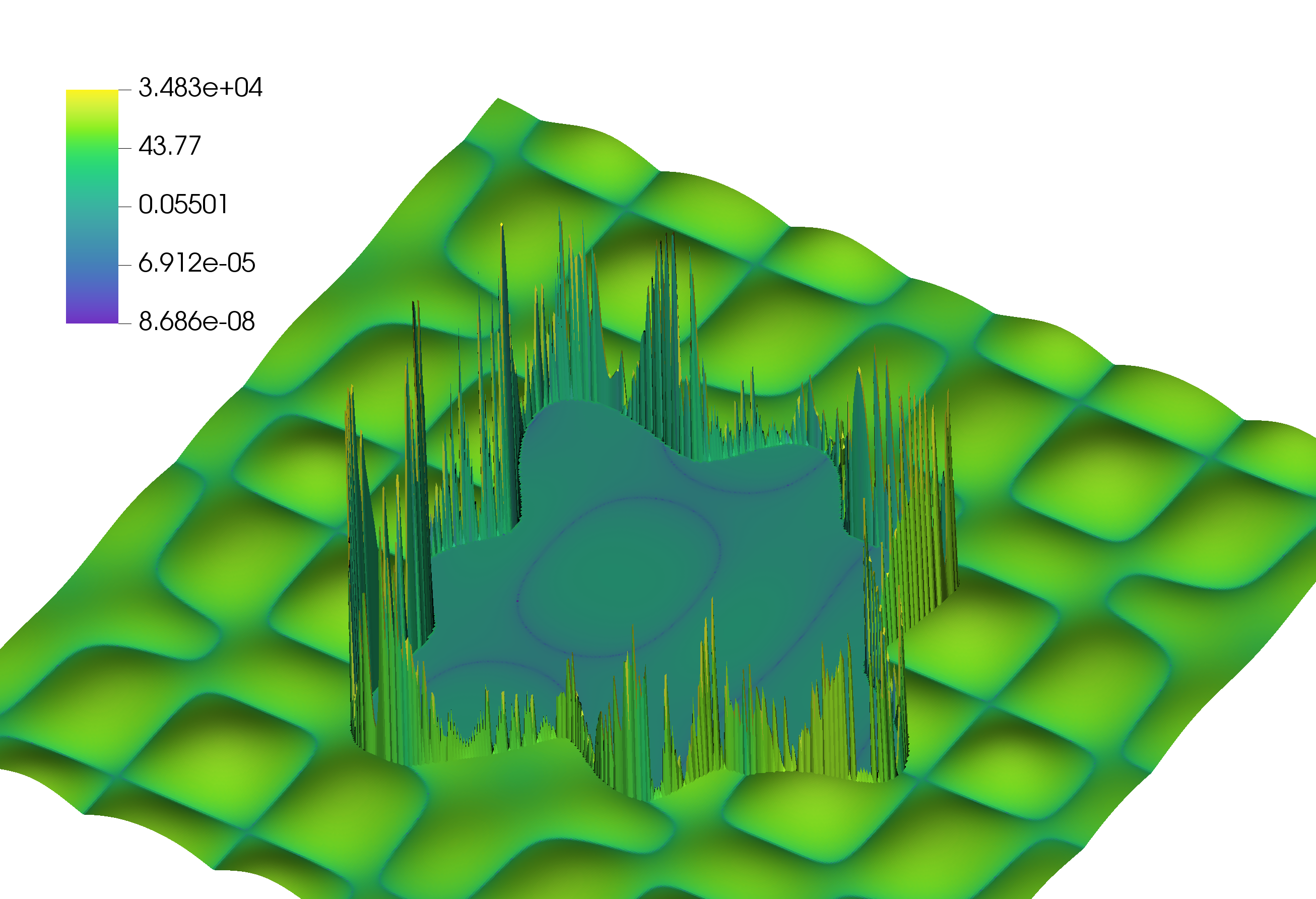}}
    \subfloat[]{
    \includegraphics[width=70mm]{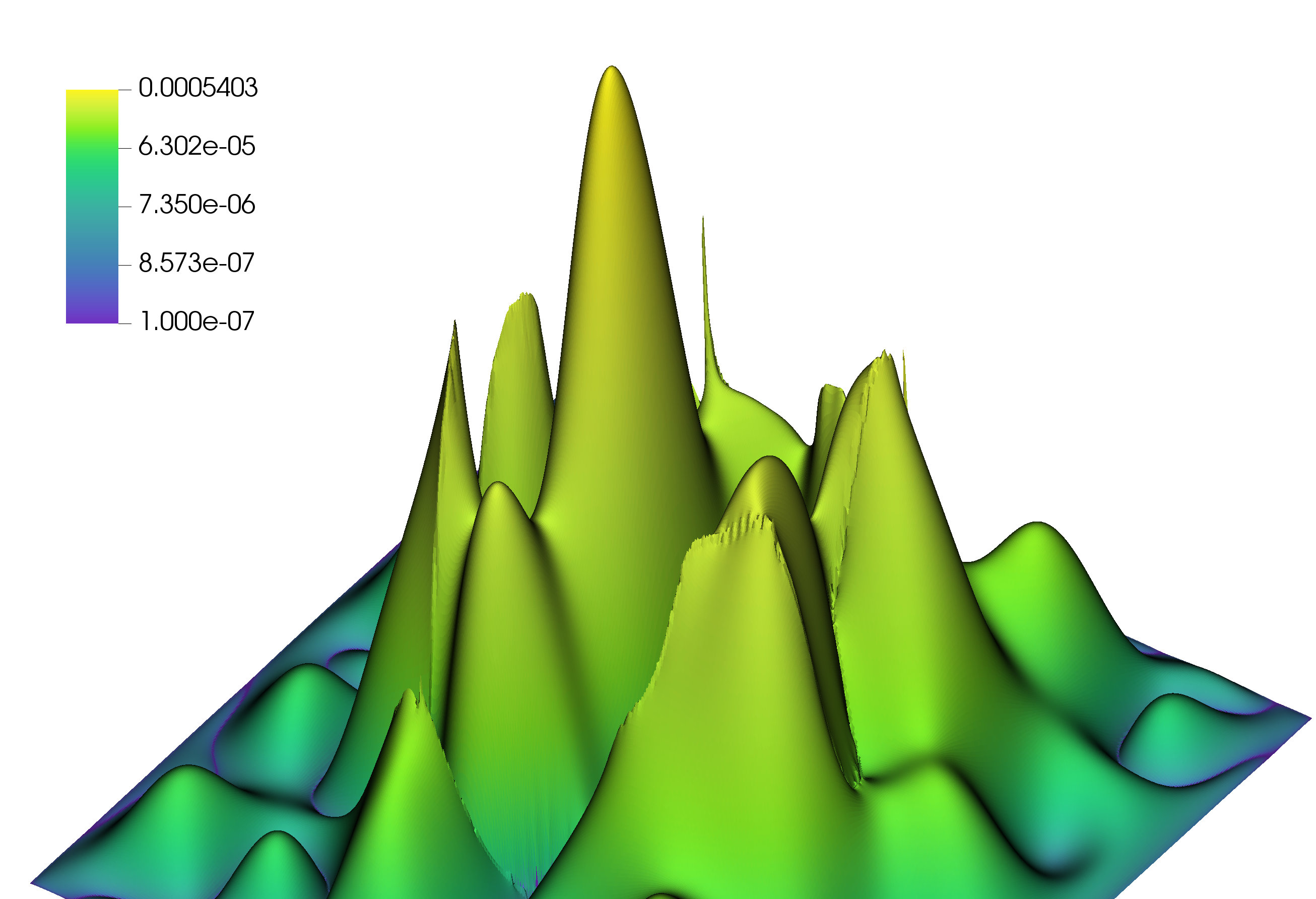}
    }
    \hspace{0 mm}
    \subfloat[]{
    \includegraphics[width=70mm]{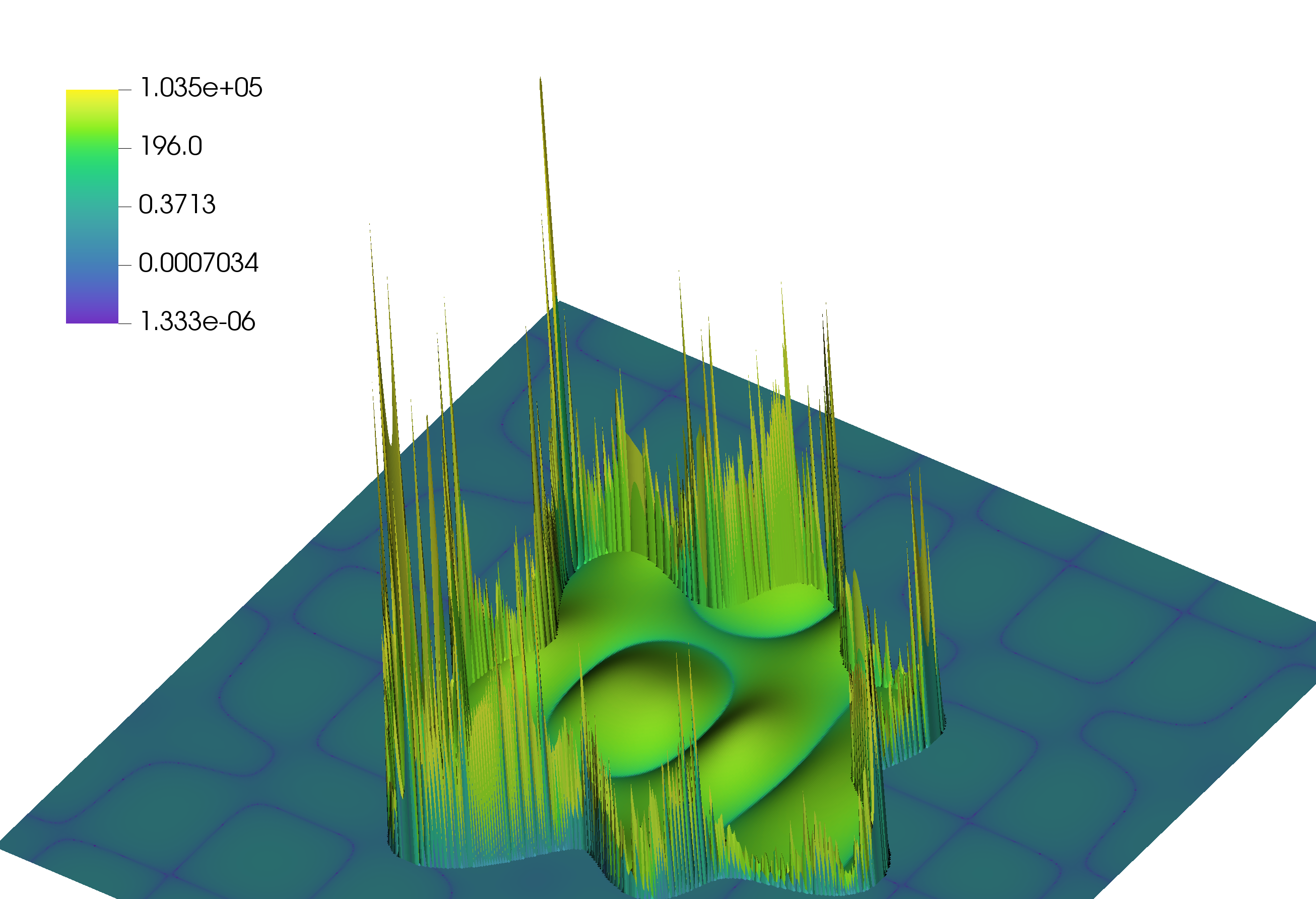}}
    \subfloat[]{
    \includegraphics[width=70mm]{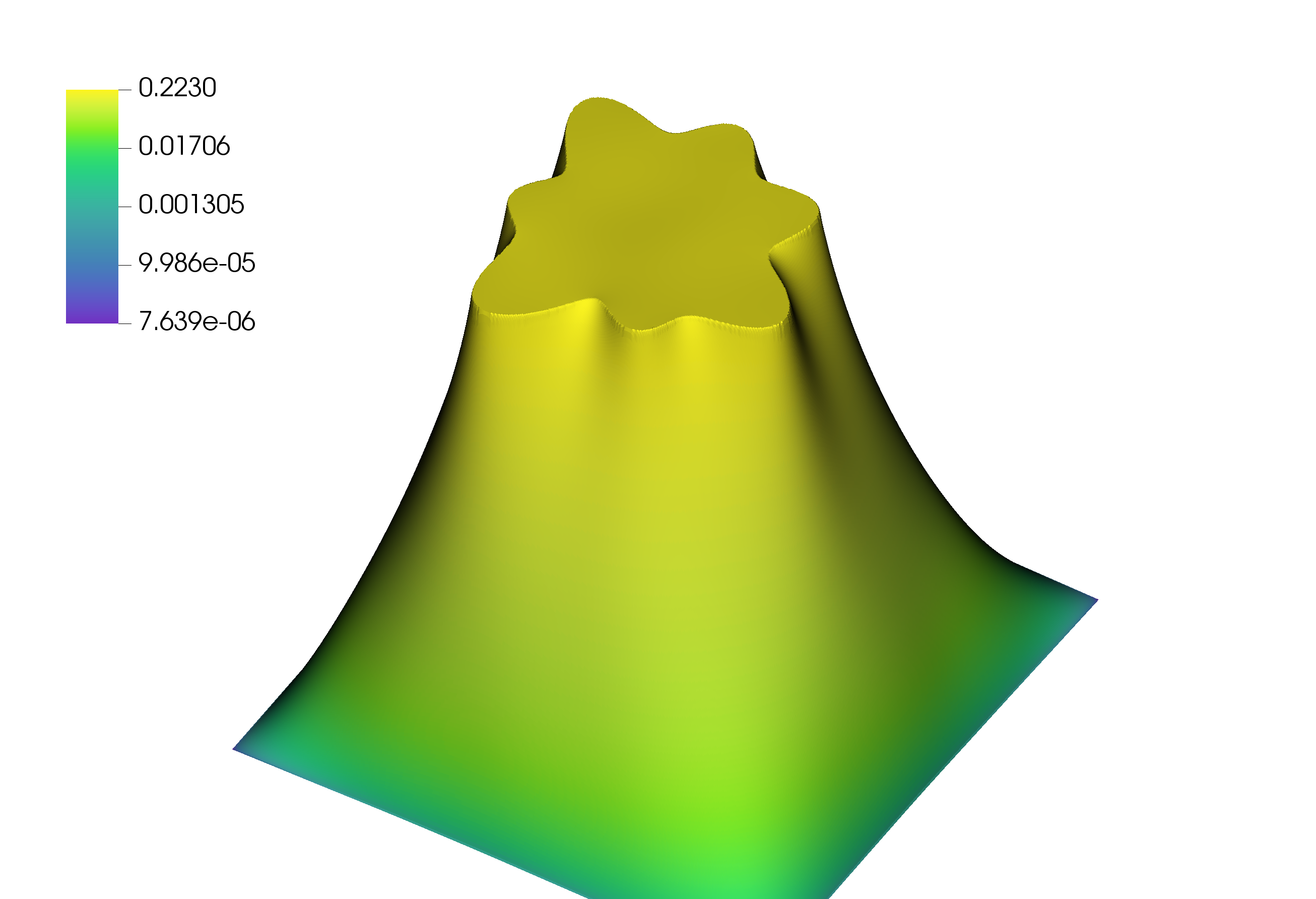}
    }
    \caption{Plots of the absolute value of truncation error (left column) and absolute value of solution error (right column) for $P=2$, $n=512$, on a log scale. Figures show: (a) and (b)  cosine geometry with $\frac{\beta^-}{\beta^+} = 10^{-4}$; (c) and (d) annulus geometry with $\frac{\beta^-}{\beta^+} = 10^{4}$; and (e) and (f) annulus geometry with $\frac{\beta^-}{\beta^+} = 10^{-4}$. Note that the larger, rough trunction error near the interface becomes smoother, rapidly decaying error in the solution (except (f), see text for discussion).}
    \label{fig:coefTestPlot}
\end{figure}

\subsection{Discontinuous Diffusion Coefficient}
Next we test the robustness of our scheme on problems with large jumps in the diffusion coefficient, as is in common in the literature (see \cite{Crockett}, \cite{gibou}, \cite{chen}, \cite{IIM}). We set the linear term coefficient $\alpha^{\pm} = 0$, and let the diffusion coefficients be constant, varying the ratio $\frac{\beta^-}{\beta^+}$ from $10^4$ to $10^{-4}$. Specifically, we fix $\beta^+ = 1$ and vary $\beta^-$ from $1$ to $10^4$, and vice versa. The manufactured solution is the same as in the truncation and solution error tests. We are interested in studying the relationship between the solution error and the ratio of diffusion coefficients, so we fix the grid spacing at $h = 256^{-1}$. However, we introduce two geometries that expose different error characteristics.

The first interface geometry is the zero level set of the function
\begin{align}
    \psi(x,y) = \frac{1}{4} \cos(\pi y) + x + \frac{\pi}{100} \ ,
\end{align}
which is simply a cosine in the xy plane. We impose periodic boundary conditions in the $y$ direction and Dirichlet boundary conditions in the $x$ direction. The $+$ region is the to the right of the interface for this geometry. We impose Dirichlet boundary conditions by filling layers of ghost cells with exact solution values. In this case both phases are tied to boundary conditions, so we do not expect any significant difference between large and small value of $\frac{\beta^-}{\beta^+}$. For the second test, we use the ``annulus" interface geometry given in section 3.1 of \cite{gibou}. This interface shape can be seen in Figures \ref{fig:coefTestPlot}(c) and (e). In the former case, the $+$ (interior) phase has a relatively large coefficient ratio ($10^4$), while in the latter it is the inverse ($10^{-4}$). Given that in these cases the $-$ phase has domain boundary conditions while the $+$ phase does not, we expect to see different error behavior as $\frac{\beta^-}{\beta^+}$ varies.

For the cosine geometry, we see that accuracy is mostly unaffected by changes in the diffusion coefficient ratio (Figure \ref{fig:coefTestError}(a)). In Figure \ref{fig:coefTestPlot}(a) and (b), the truncation error again is larger at the interface but the solution error is smooth, similar to the previous tests. We observed no impact of the conditioning of the discretization matrix on the solution accuracy as the ratio of coefficients varies over 8 orders of magnitude.

For the second test, instead we see that the error is about 4 orders of magnitude higher when we have the diffusion coefficient on the interior of the domain is much larger than the diffusion coefficient on the exterior of the domain ($\beta^- / \beta^+ = 10^{-4}$). This result is consistent with results reported in Figure 2 of \cite{Crockett},  as well as Figure 17b of \cite{chen}. Through potential theoretic arguments, we believe this is a result of solution errors in the interior region not being ``tied down'' to any domain boundary condition, as in the cosine test. Because of the gradient jump condition \eqref{eq:dujump}, any gradient \emph{errors} in the interior are multiplied by $\beta^+ / \beta^- = 10^4$ in their contribution to the exterior domain gradients at the interface, forcing the interior solution there to ``drift'' in proportion. However, with this scaling the convergence rates are still retained, but with an error constant reflecting this ratio in diffusion coefficients.

\FloatBarrier

\subsection{Discontinuous Solution}
\begin{figure}[h]
\centering
\subfloat[]{
\begin{tikzpicture}[scale=.9]
    \begin{axis}[
    ymode = log,
    xmode = log,
    title={Solution Error, $L^{\infty}$ Norm},
    xlabel={$\frac{s^-}{s^+}$},
    ylabel={Error},
    xmin=1e-4, xmax=1e4,
    ymin=1e-6, ymax=1e4,
    xtick={1e-4, 1e-2,  1e0,  1e2,  1e4},
    xticklabels={1e-4, 1e-2,  1e0,  1e2,  1e4},
    ytick={1e-6, 1e-4, 1e-2, 1e-0,1e2,1e4},
    legend pos=south west,
    ymajorgrids=true,
    grid style=dashed,
]

\addplot[
    color=cyan,
    mark=square,
    ]
    coordinates {
    (1e-4, 1.728e+02) (1e-3,1.727e+01 ) (1e-2,1.721e+00 )(1e-1, 1.659e-01) (1e0,1.042e-02 )(1e1,8.284e-02 )(1e2,8.275e-01 )(1e3,8.274e+00 )(1e4,8.274e+01 )
    };
    \addlegendentry{$P=2$}

\addplot[
    color=red,
    mark=square,
    ]
    coordinates {
    (1e-4, 4.155e+00) (1e-3,4.152e-01 ) (1e-2,4.121e-02 )(1e-1,3.809e-03) (1e0,3.673e-04 )(1e1,3.625e-03  )(1e2,3.652e-02 )(1e3,3.655e-01 )(1e4,3.655e+00 )
    };
    \addlegendentry{$P=4$}
\addplot[
    color=blue,
    mark=square,
    ]
    coordinates {
    (1e-4, 1.662e-01) (1e-3,1.661e-02 ) (1e-2,1.648e-03 )(1e-1,1.516e-04) (1e0,1.570e-05 )(1e1,1.577e-04  )(1e2,1.601e-03 )(1e3, 1.603e-02  )(1e4,1.603e-01  )
    };
    \addlegendentry{$P=6$}

\end{axis}
\end{tikzpicture}
}
\subfloat[]{
\begin{tikzpicture}[scale=.9]
    \begin{axis}[
    ymode = log,
    xmode = log,
    title={Solution Error, $L^{\infty}$ Norm},
    xlabel={$\frac{s^-}{s^+}$},
    ylabel={Error},
    xmin=1e-4, xmax=1e4,
    ymin=1e-6, ymax=1e4,
    xtick={1e-4, 1e-2,  1e0,  1e2,  1e4},
    xticklabels={1e-4, 1e-2,  1e0,  1e2,  1e4},
    ytick={1e-6, 1e-4, 1e-2, 1e-0,1e2,1e4},
    legend pos=south west,
    ymajorgrids=true,
    grid style=dashed,
]

\addplot[
    color=cyan,
    mark=square,
    ]
    coordinates {
    (1e-4, 1.585e+02) (1e-3,1.585e+01) (1e-2,1.585e+00  )(1e-1, 1.586e-01) (1e0,1.591e-02 )(1e1,4.722e-02 )(1e2,4.192e-01 )(1e3,4.139e+00 )(1e4,4.134e+01 )
    };
    \addlegendentry{$P=2$}

\addplot[
    color=red,
    mark=square,
    ]
    coordinates {
    (1e-4, 4.023e+00) (1e-3,4.022e-01 ) (1e-2,4.011e-02  )(1e-1,3.908e-03) (1e0,3.576e-04 )(1e1,2.147e-03  )(1e2,2.022e-02  )(1e3,2.011e-01 )(1e4,2.010e+00 )
    };
    \addlegendentry{$P=4$}
\addplot[
    color=blue,
    mark=square,
    ]
    coordinates {
    (1e-4, 1.896e-01) (1e-3,1.896e-02) (1e-2,1.898e-03 )(1e-1,1.915e-04 ) (1e0,2.079e-05 )(1e1,1.230e-04  )(1e2,1.178e-03)(1e3, 1.173e-02 )(1e4,1.173e-01  )
    };
    \addlegendentry{$P=6$}

\end{axis}
\end{tikzpicture}
}




\caption{Comparison of solution maximum errors for two geometries (a) cosine, and (b) annulus, and over a range of solution scaling ratios, fixing grid size at $n=128$ cells and varying the order of accuracy $P$ (note that the $L^1$ error is nearly identical).}

\label{fig:solTest}
\end{figure}
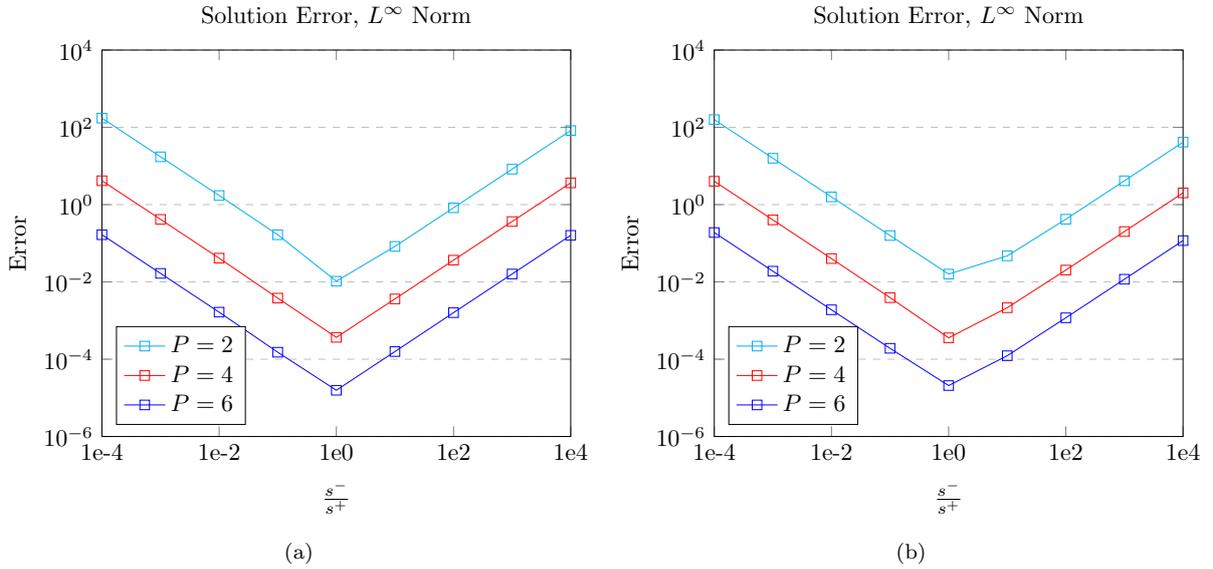

We perform a similar test with a solution that has large jumps at the interface. $\beta^{\pm}$ and $u^{\pm}$ are the same as in the solution and truncation error test, and we again set $\alpha^{\pm}=0$. The $u^{\pm}$ fields are multiplied by scaling factors $s^{\pm}$ to create large jumps in the solution, and we use the same two geometries as the discontinuous coefficients test. We observe $L^1$ and $L^{\infty}$ errors that are proportional to the larger of the two scaling coefficients (see Figure \ref{fig:solTest}). This scaling does not magnify the error because it appears as a large discontinuity in the source term, as well as in the jump conditions, which both contribute only to the right hand side of the linear system. This test highlights the importance of the having two separate degrees of freedom in each cut cell, from which we are able to accurately reproduce a solution and gradients which jump by up to 4 orders of magnitude across the interface.

\FloatBarrier

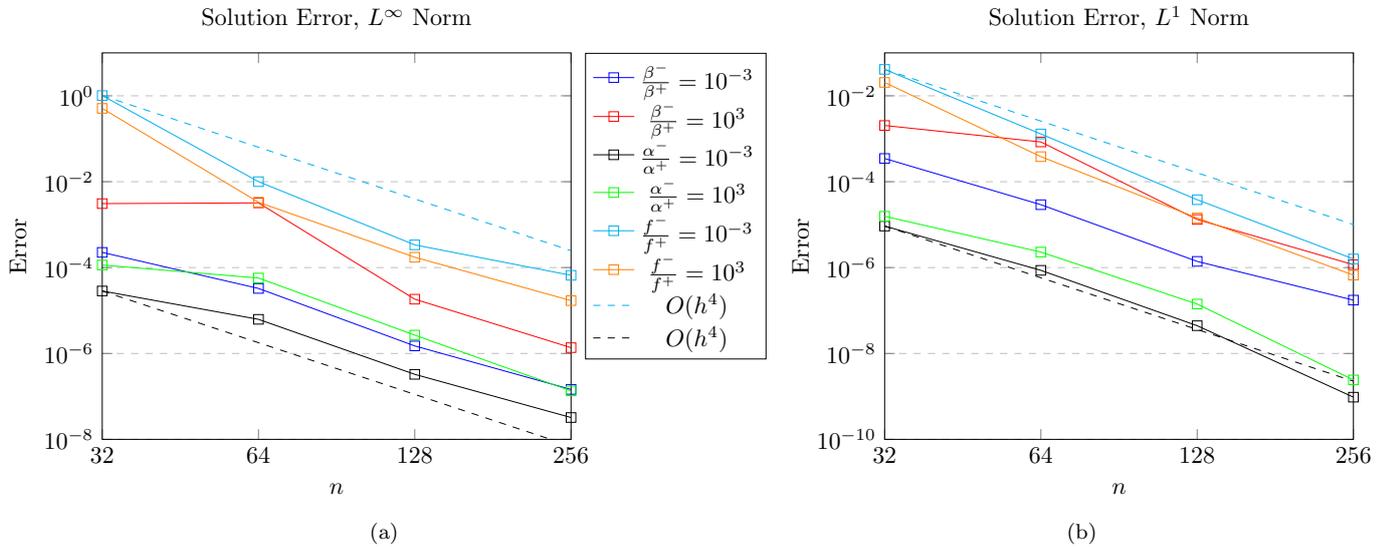
\begin{figure}
\centering
\subfloat[]{
\begin{tikzpicture}[scale=.9]
\begin{axis}[
    xmode = log,
    ymode = log,
    title={Solution Error, $L^{\infty}$ Norm},
    xlabel={$n$},
    ylabel={Error},
    xmin=32, xmax=256,
    ymin=1e-8, ymax=1e1,
    xtick={32,64,128,256},
    xticklabels={32,64,128,256},
    ytick={1e-8,1e-6,1e-4,1e-2,1e0},
    legend pos= outer north east,
    ymajorgrids=true,
    grid style=dashed,
]

\addplot[
    color=blue,
    mark=square,
    ]
    coordinates {
    (32,2.259e-04)(64,3.269e-05) (128,1.498e-06)(256,1.452e-07)
    };
    \addlegendentry{$\frac{\beta^-}{\beta^+} = 10^{-3}$ }
\addplot[
    color=red,
    mark=square,
    ]
    coordinates {
    (32,3.090e-03)(64,3.211e-03) (128,1.837e-05)(256,1.366e-06)
    };
    \addlegendentry{$\frac{\beta^-}{\beta^+} = 10^{3}$ }

\addplot[
    color=black,
    mark=square,
    ]
    coordinates {
    (32,2.851e-05)(64,6.234e-06) (128,3.265e-07 )(256,3.224e-08)
    };
    \addlegendentry{$\frac{\alpha^-}{\alpha^+} = 10^{-3}$ }
\addplot[
    color=green,
    mark=square,
    ]
    coordinates {
    (32,1.148e-04)(64,5.717e-05) (128,2.675e-06)(256,1.345e-07)
    };
    \addlegendentry{$\frac{\alpha^-}{\alpha^+} = 10^{3}$ }

\addplot[
    color=cyan,
    mark=square,
    ]
    coordinates {
    (32,1.015e+00)(64,1.003e-02 ) (128,3.380e-04)(256,6.607e-05)
    };
    \addlegendentry{$\frac{f^-}{f^+} = 10^{-3}$ }
\addplot[
    color=orange,
    mark=square,
    ]
    coordinates {
    (32,5.102e-01)(64,3.326e-03) (128,1.741e-04)(256,1.701e-05)
    };
    \addlegendentry{$\frac{f^-}{f^+} = 10^{3}$}

\addplot[
    color=cyan,
    dashed,
    ]
    coordinates {
    (32,1.015e+00)(256,2.4780e-04)
    };
    \addlegendentry{$O(h^4)$}
\addplot[
    color=black,
    dashed,
    ]
    coordinates {
    (32,2.851e-05)(256,6.9604e-09)
    };
    \addlegendentry{$O(h^4)$}

\end{axis}
\end{tikzpicture}
}
\subfloat[]{
\begin{tikzpicture}[scale=.9]
\begin{axis}[
    xmode = log,
    ymode = log,
    title={Solution Error, $L^{1}$ Norm},
    xlabel={$n$},
    ylabel={Error},
    xmin=32, xmax=256,
    ymin=1e-10, ymax=1e-1,
    xtick={32,64,128,256},
    xticklabels={32,64,128,256},
    ytick={1e-10,1e-8,1e-6,1e-4,1e-2},
    ymajorgrids=true,
    grid style=dashed,
]

\addplot[
    color=blue,
    mark=square,
    ]
    coordinates {
    (32,3.475e-04)(64,2.916e-05 ) (128,1.394e-06)(256,1.747e-07)
    };
\addplot[
    color=red,
    mark=square,
    ]
    coordinates {
    (32,2.038e-03)(64,8.436e-04) (128,1.319e-05)(256,1.166e-06)
    };

\addplot[
    color=black,
    mark=square,
    ]
    coordinates {
    (32,9.265e-06)(64,8.626e-07) (128,4.434e-08 )(256,9.604e-10)
    };
\addplot[
    color=green,
    mark=square,
    ]
    coordinates {
    (32,1.567e-05)(64,2.298e-06) (128,1.414e-07)(256,2.423e-09 )
    };

\addplot[
    color=cyan,
    mark=square,
    ]
    coordinates {
    (32,4.125e-02)(64,1.294e-03  ) (128,3.808e-05)(256,1.589e-06)
    };
\addplot[
    color=orange,
    mark=square,
    ]
    coordinates {
    (32,2.050e-02)(64,3.832e-04) (128,1.432e-05)(256,6.676e-07)
    };

\addplot[
    color=cyan,
    dashed,
    ]
    coordinates {
    (32,4.125e-02)(256,0.0000100708)
    };
\addplot[
    color=black,
    dashed,
    ]
    coordinates {
    (32,9.265e-06)(256,2.2620e-09)
    };

\end{axis}
\end{tikzpicture}
}
\caption{$L^{\infty}$ and $L^1$ norm of errors for Richardson convergence test for a variety of scaling factor ratios in linear term coefficient, diffusion coefficient, and source term. We fix $P=4$ and observe fourth order convergence in the $L^1$ norm and nearly fourth order convergence in the max norm. }
\label{fig:geomError}
\end{figure}

\subsection{Imposing Homogeneous Jump Conditions}
Lastly, we test the ability of the method to impose jump conditions as a constraint. We let $\alpha^{\pm}, \beta^{\pm}$ be the same as in the truncation and solution error test, and use the manufactured solution $u^{\pm}$ from that test as the source term $f^{\pm}$.  We impose homogeneous (zero) jump conditions and focus on the annulus geometry. For this test we fix $P=4$ and test the scheme with a variety of coefficient and source term scalings. The error is measured by using the $n=512$ numerical solution as the exact solution, with the results in Figure \ref{fig:geomError}. We observe roughly fourth order convergence for all tests, although there is more variation in convergence compared to the manufactured solution tests.
\FloatBarrier

\section{Conclusion}

We have developed a finite volume method for the variable coefficient elliptic interface problem and demonstrated up to sixth order accurate on a variety of test problems. In developing this method we gave a general truncation error analysis that justifies the use of stencils based on least-squares interpolation. Our stencils are derived from cell-centered Taylor polynomials which are implicitly defined in terms of local values of the solution and interface jumps or boundary conditions, where appropriate. To enforce conservation, we choose a single flux on each face which is an average of the flux calculated in neighboring cells from the respective Taylor polynomials. In cells away from the interface, we take advantage of standard finite volume symmetries to build stencils with a minimal footprint. 

Future research will involve 1) building an efficient geometric multigrid solver along the lines of \cite{Devendran2015AHM}, 2) extending our method to three dimensions, and 3) incorporating adaptive mesh refinement. The prior method presented in \cite{Dharshi} accomplishes these for Poisson's equation in more complex geometries, including boundaries with kinks, using ``smoothed'' constructive solid geometry capabilities of the Chombo software library \cite{Chombo}. This would enable this method to be used in discretizations for large scale science applications. Further exploration is also needed of the theory of undetermined stencil systems (using the moment matrix transpose, $\*M^T$). This paper has shown that building a stencil of a given order and truncation error still allows infinitely many valid stencils; this fact could be exploited to promote sparsity or alter the conditioning or stability of the operator, and we are drafting a paper with analysis that may provide specific algorithmic guidance.

\section*{Funding:}
This work is supported by the U.S. Department of Energy, Office of Science, Advanced Scientific Computing Research, Base Math Program, under contract number DE-AC02-05CH11231.

\appendix
\section{Geometry Generation Algorithm} 
We specify the interface as a zero level set of an implicit function $\psi(\*x)$. Therefore, for sufficiently smooth $\phi$ compared to the grid resolution, we assume we can identify cut cells by evaluating $\psi$ at the four corners of each cell. If any of these values have different signs, the cell is tagged as a cut cell. If the interface intersects one face of the cell multiple times, or there are more than two faces intersected by the interface, we consider the geometry to be under-resolved and could refine the mesh or adjust the boundary without inducing significant errors. Given these assumptions, when the interface intersects a cell it creates a region which is bounded on one side by the interface and on two or three sides by the edges of a square cell (see Figure \ref{fig:cutCell}). Volume moments \eqref{eq:vMom} are defined as integrals over this region. Area moments, defined in \eqref{eq:EBMom} and \eqref{eq:nMom}, are integrals over the portion of the EB that intersects the cut cell. We compute these integrals by approximating the interface with piecewise line segments, and then apply a formula for the integral of monomials along line segments. The vertices of the line segments are roots of $\psi$, which we find using a simple root finder such as the secant method. By refining this interface iteratively into $2^n$ line segments, we can calculate a convergent sequence $m_n$ of moment approximations that stops when $\abs{m_{n+1}- m_n}$ reaches machine precision. The convergence of this sequence is accelerated using Richardson extrapolation, which in this case is often referred to as Romberg integration.  

A formula for the integral of $x^py^q$ over an arbitrary polygon can be derived from Green's theorem:
\begin{align}
    \int_P \frac{\partial f(x,y)}{\partial x} dV = \int_C f(x,y) \hat{n}_x dA = \int_C f(x,y) dy \ ,
\end{align}
where $C$ is the boundary of the polygon $P$. Let $f = \frac{x^{p+1}}{p+1}y^q$, giving us:
\begin{align}
    \int_P x^py^q dV = \int_C \frac{x^{p+1}}{p+1}y^q dy \ . \label{eq:int}
\end{align}
We parameterize each edge segment by:
\begin{align}
    x(t) &= (x_{k+1}-x_k)t + x_k = \Delta x t+x_k\\
    y(t) &= (y_{k+1}-y_k)t + y_k = \Delta y t + y_k \ , \label{eq:linearParam}
\end{align}
where $t$ goes from $0$ to $1$ and $(x_k,y_k)$ are the ordered vertices of the polygon. The formula for the right hand side of \eqref{eq:int} along a single line segment $C_k(t)$ is obtained by a binomial expansion:
\begin{align}
    \int_{C_k} \frac{x^{p+1}}{p+1}y^q dy &=  \frac{1}{p+1}\int_0^1 (\Delta x t + x_k)^{p+1} (\Delta yt + y_k)^q \Delta y dt \\
    &= \sum_{i=0}^{p+1} \sum_{j=0}^{q} \frac{{p+1 \choose i} {q \choose j}}{(p+1)(p+2+q-i-j)} (x_k^iy_k^j)(\Delta x)^{p+1-i}(\Delta y)^{q+1-j} \ .
\end{align}
We follow a similar procedure for area integrals:
\begin{align}
    \int_C x^p y^q &= \int_0^1 x(t)^p y(t)^q dC \\
    &= \sum_{i=0}^{p} \sum_{j=0}^{q} \frac{{p \choose i} {q \choose j}}{(p+1+q-i-j)} (x_k^iy_k^j)(\Delta x)^{p-i}(\Delta y)^{q-j} \Delta C\ ,  \label{eq:area}
\end{align}
and to calculate area integrals times unit normals we multiply equation \eqref{eq:area} by $n_x = \dfrac{\Delta y}{\Delta C}$ or $n_y = \dfrac{-\Delta x}{\Delta C}$, where $\Delta C = \sqrt{\Delta x^2 + \Delta y^2}$. (The tangent vector is rotated 90 degrees clockwise). The integrals of interest are obtained by adding the integrals along all line segments of the polygon in the case of volume moments, or just along the interface in the case of area moments.


\newpage    

 \bibliographystyle{elsarticle-num-names} 
 \bibliography{cas-refs}





\end{document}

%% file: abbrev.tex
\def\*#1{\mathbf{#1}}
\def\avg#1{\langle #1 \rangle}

\newcommand{\lrp}[1]{\left( #1 \right)}
\newcommand{\lrb}[1]{\left[ #1 \right ]}
\newcommand{\IP}[1]{\left \langle #1 \right \rangle}
\newcommand{\IPi}[1]{\left \langle #1 \right \rangle_{\*i}}

\newcommand{\Vi}{V_{\*i}}
\newcommand{\Vj}{V_{\*j}}

\newcommand{\Vpi}{V_{p,\*i}}
\renewcommand{\set}[1]{\left\{#1\right\}}

\renewcommand{\l}{\bm{\ell}}
\newcommand{\abs}[1]{\left|#1\right|}
\newcommand{\Norm}[1]{\left\Vert#1\right\Vert}

%% file: main.bbl
\begin{thebibliography}{23}
\expandafter\ifx\csname natexlab\endcsname\relax\def\natexlab#1{#1}\fi
\providecommand{\url}[1]{\texttt{#1}}
\providecommand{\href}[2]{#2}
\providecommand{\path}[1]{#1}
\providecommand{\DOIprefix}{doi:}
\providecommand{\ArXivprefix}{arXiv:}
\providecommand{\URLprefix}{URL: }
\providecommand{\Pubmedprefix}{pmid:}
\providecommand{\doi}[1]{\href{http://dx.doi.org/#1}{\path{#1}}}
\providecommand{\Pubmed}[1]{\href{pmid:#1}{\path{#1}}}
\providecommand{\bibinfo}[2]{#2}
\ifx\xfnm\relax \def\xfnm[#1]{\unskip,\space#1}\fi
\bibitem[{Li(2003)}]{LiReview}
\bibinfo{author}{Z.~Li},
\newblock \bibinfo{title}{An overview of the immersed interface method and its
  applications},
\newblock \bibinfo{journal}{Taiwanese Journal of Mathematics}
  \bibinfo{volume}{7} (\bibinfo{year}{2003}) \bibinfo{pages}{1--49}.
\bibitem[{Gibou et~al.(2013)Gibou, Min, and Fedkiw}]{gibouReview}
\bibinfo{author}{F.~Gibou}, \bibinfo{author}{C.~Min},
  \bibinfo{author}{R.~Fedkiw},
\newblock \bibinfo{title}{High resolution sharp computational methods for
  elliptic and parabolic problems in complex geometries},
\newblock \bibinfo{journal}{J. Sci. Comput.} \bibinfo{volume}{54}
  (\bibinfo{year}{2013}) \bibinfo{pages}{369–413}.
\bibitem[{Babuska(1970)}]{Babuska1970TheFE}
\bibinfo{author}{I.~Babuska},
\newblock \bibinfo{title}{The finite element method for elliptic equations with
  discontinuous coefficients},
\newblock \bibinfo{journal}{Computing} \bibinfo{volume}{5}
  (\bibinfo{year}{1970}) \bibinfo{pages}{207--213}.
\bibitem[{Li(1998)}]{IIMFE}
\bibinfo{author}{Z.~Li},
\newblock \bibinfo{title}{The immersed interface method using a finite element
  formulation},
\newblock \bibinfo{journal}{Applied Numerical Mathematics} \bibinfo{volume}{27}
  (\bibinfo{year}{1998}) \bibinfo{pages}{253--267}.
\bibitem[{Leveque and Li(1994)}]{IIM}
\bibinfo{author}{R.~J. Leveque}, \bibinfo{author}{Z.~Li},
\newblock \bibinfo{title}{The immersed interface method for elliptic equations
  with discontinuous coefficients and singular sources},
\newblock \bibinfo{journal}{SIAM Journal on Numerical Analysis}
  \bibinfo{volume}{31} (\bibinfo{year}{1994}) \bibinfo{pages}{1019--1044}.
\bibitem[{Fedkiw et~al.(1999)Fedkiw, Aslam, Merriman, and Osher}]{GFM}
\bibinfo{author}{R.~P. Fedkiw}, \bibinfo{author}{T.~Aslam},
  \bibinfo{author}{B.~Merriman}, \bibinfo{author}{S.~Osher},
\newblock \bibinfo{title}{A non-oscillatory eulerian approach to interfaces in
  multimaterial flows (the ghost fluid method)},
\newblock \bibinfo{journal}{Journal of Computational Physics}
  \bibinfo{volume}{152} (\bibinfo{year}{1999}) \bibinfo{pages}{457--492}.
\bibitem[{Johansen and Colella(1998)}]{Hans_thesis}
\bibinfo{author}{H.~Johansen}, \bibinfo{author}{P.~Colella},
\newblock \bibinfo{title}{A cartesian grid embedded boundary method for
  poisson's equation on irregular domains},
\newblock \bibinfo{journal}{Journal of Computational Physics}
  \bibinfo{volume}{147} (\bibinfo{year}{1998}) \bibinfo{pages}{60--85}.
\bibitem[{Crockett et~al.(2011)Crockett, Colella, and Graves}]{Crockett}
\bibinfo{author}{R.~Crockett}, \bibinfo{author}{P.~Colella},
  \bibinfo{author}{D.~Graves},
\newblock \bibinfo{title}{A cartesian grid embedded boundary method for solving
  the poisson and heat equations with discontinuous coefficients in three
  dimensions},
\newblock \bibinfo{journal}{Journal of Computational Physics}
  \bibinfo{volume}{230} (\bibinfo{year}{2011}) \bibinfo{pages}{2451--2469}.
\bibitem[{SCH(2006)}]{SCHWARTZ2006531}
\bibinfo{title}{A cartesian grid embedded boundary method for the heat equation
  and poisson’s equation in three dimensions},
\newblock \bibinfo{journal}{Journal of Computational Physics}
  \bibinfo{volume}{211} (\bibinfo{year}{2006}) \bibinfo{pages}{531--550}.
\bibitem[{Colella(2016)}]{HOFV}
\bibinfo{author}{P.~Colella},
\newblock \bibinfo{title}{High-order finite-volume methods on
  locally-structured grids},
\newblock \bibinfo{journal}{Discrete and Continuous Dynamical Systems}
  \bibinfo{volume}{36} (\bibinfo{year}{2016}) \bibinfo{pages}{4247--4270}.
\bibitem[{Chen and Strain(2008)}]{chen}
\bibinfo{author}{T.~Chen}, \bibinfo{author}{J.~Strain},
\newblock \bibinfo{title}{Piecewise-polynomial discretization and
  krylov-accelerated multigrid for elliptic interface problems},
\newblock \bibinfo{journal}{Journal of Computational Physics}
  \bibinfo{volume}{227} (\bibinfo{year}{2008}) \bibinfo{pages}{7503--7542}.
\bibitem[{Devendran et~al.(2017)Devendran, Graves, Johansen, and
  Ligocki}]{Dharshi}
\bibinfo{author}{D.~Devendran}, \bibinfo{author}{D.~Graves},
  \bibinfo{author}{H.~Johansen}, \bibinfo{author}{T.~Ligocki},
\newblock \bibinfo{title}{{A fourth-order Cartesian grid embedded boundary
  method for Poisson's equation}},
\newblock \bibinfo{journal}{Communications in Applied Mathematics and
  Computational Science} \bibinfo{volume}{12} (\bibinfo{year}{2017})
  \bibinfo{pages}{51 -- 79}.
\bibitem[{Zhang et~al.(2012)Zhang, Johansen, and Colella}]{Zhang}
\bibinfo{author}{Q.~Zhang}, \bibinfo{author}{H.~Johansen},
  \bibinfo{author}{P.~Colella},
\newblock \bibinfo{title}{A fourth-order accurate finite-volume method with
  structured adaptive mesh refinement for solving the advection-diffusion
  equation},
\newblock \bibinfo{journal}{SIAM Journal on Scientific Computing}
  \bibinfo{volume}{34} (\bibinfo{year}{2012}) \bibinfo{pages}{B179--B201}.
\bibitem[{LeVeque(2007)}]{leveque}
\bibinfo{author}{R.~J. LeVeque}, \bibinfo{title}{Finite Difference Methods for
  Ordinary and Partial Differential Equations}, \bibinfo{publisher}{Society for
  Industrial and Applied Mathematics}, \bibinfo{year}{2007}.
\bibitem[{Overton-Katz et~al.(2022)Overton-Katz, Gao, Guzik, Antepara, Graves,
  and Johansen}]{nate}
\bibinfo{author}{N.~Overton-Katz}, \bibinfo{author}{X.~Gao},
  \bibinfo{author}{S.~Guzik}, \bibinfo{author}{O.~Antepara},
  \bibinfo{author}{D.~T. Graves}, \bibinfo{author}{H.~Johansen},
\newblock \bibinfo{title}{A fourth-order embedded boundary finite volume method
  for the unsteady stokes equations with complex geometries},
\newblock \bibinfo{journal}{arXiv}  (\bibinfo{year}{2022}).
\bibitem[{Barker et~al.(2001)Barker, Blackford, Dongarra, Croz, Hammarling,
  Marinova, Waśniewski, and Yalamov}]{LAPACK}
\bibinfo{author}{V.~A. Barker}, \bibinfo{author}{L.~S. Blackford},
  \bibinfo{author}{J.~Dongarra}, \bibinfo{author}{J.~D. Croz},
  \bibinfo{author}{S.~Hammarling}, \bibinfo{author}{M.~Marinova},
  \bibinfo{author}{J.~Waśniewski}, \bibinfo{author}{P.~Yalamov},
  \bibinfo{title}{LAPACK95 Users' Guide}, \bibinfo{publisher}{Society for
  Industrial and Applied Mathematics}, \bibinfo{year}{2001}.
\bibitem[{Balay et~al.(2022)Balay, Abhyankar, Adams, Benson, Brown, Brune,
  Buschelman, Constantinescu, Dalcin, Dener, Eijkhout, Gropp, Hapla, Isaac,
  Jolivet, Karpeev, Kaushik, Knepley, Kong, Kruger, May, McInnes, Mills,
  Mitchell, Munson, Roman, Rupp, Sanan, Sarich, Smith, Zampini, Zhang, Zhang,
  and Zhang}]{petsc-user-ref}
\bibinfo{author}{S.~Balay}, \bibinfo{author}{S.~Abhyankar},
  \bibinfo{author}{M.~F. Adams}, \bibinfo{author}{S.~Benson},
  \bibinfo{author}{J.~Brown}, \bibinfo{author}{P.~Brune},
  \bibinfo{author}{K.~Buschelman}, \bibinfo{author}{E.~Constantinescu},
  \bibinfo{author}{L.~Dalcin}, \bibinfo{author}{A.~Dener},
  \bibinfo{author}{V.~Eijkhout}, \bibinfo{author}{W.~D. Gropp},
  \bibinfo{author}{V.~Hapla}, \bibinfo{author}{T.~Isaac},
  \bibinfo{author}{P.~Jolivet}, \bibinfo{author}{D.~Karpeev},
  \bibinfo{author}{D.~Kaushik}, \bibinfo{author}{M.~G. Knepley},
  \bibinfo{author}{F.~Kong}, \bibinfo{author}{S.~Kruger},
  \bibinfo{author}{D.~A. May}, \bibinfo{author}{L.~C. McInnes},
  \bibinfo{author}{R.~T. Mills}, \bibinfo{author}{L.~Mitchell},
  \bibinfo{author}{T.~Munson}, \bibinfo{author}{J.~E. Roman},
  \bibinfo{author}{K.~Rupp}, \bibinfo{author}{P.~Sanan},
  \bibinfo{author}{J.~Sarich}, \bibinfo{author}{B.~F. Smith},
  \bibinfo{author}{S.~Zampini}, \bibinfo{author}{H.~Zhang},
  \bibinfo{author}{H.~Zhang}, \bibinfo{author}{J.~Zhang},
  \bibinfo{title}{{PETSc/TAO} Users Manual}, \bibinfo{type}{Technical Report}
  \bibinfo{number}{ANL-21/39 - Revision 3.17}, Argonne National Laboratory,
  \bibinfo{year}{2022}.
\bibitem[{Balay et~al.(1997)Balay, Gropp, McInnes, and Smith}]{petsc-efficient}
\bibinfo{author}{S.~Balay}, \bibinfo{author}{W.~D. Gropp},
  \bibinfo{author}{L.~C. McInnes}, \bibinfo{author}{B.~F. Smith},
\newblock \bibinfo{title}{Efficient management of parallelism in object
  oriented numerical software libraries},
\newblock in: \bibinfo{editor}{E.~Arge}, \bibinfo{editor}{A.~M. Bruaset},
  \bibinfo{editor}{H.~P. Langtangen} (Eds.), \bibinfo{booktitle}{Modern
  Software Tools in Scientific Computing}, \bibinfo{publisher}{Birkh{\"{a}}user
  Press}, \bibinfo{year}{1997}, pp. \bibinfo{pages}{163--202}.
\bibitem[{Li and Demmel(2003)}]{lidemmel03}
\bibinfo{author}{X.~S. Li}, \bibinfo{author}{J.~W. Demmel},
\newblock \bibinfo{title}{{SuperLU\_DIST}: A scalable distributed-memory sparse
  direct solver for unsymmetric linear systems},
\newblock \bibinfo{journal}{ACM Trans. Mathematical Software}
  \bibinfo{volume}{29} (\bibinfo{year}{2003}) \bibinfo{pages}{110--140}.
\bibitem[{Devendran et~al.(2014)Devendran, Graves, and
  Johansen}]{Devendran2015AHM}
\bibinfo{author}{D.~Devendran}, \bibinfo{author}{D.~T. Graves},
  \bibinfo{author}{H.~Johansen},
\newblock \bibinfo{title}{A hybrid multigrid algorithm for poisson’s equation
  using an adaptive, fourth order treatment of cut cells},
\newblock \bibinfo{journal}{Technical Report LBNL-1004329, LBNL}
  (\bibinfo{year}{2014}).
\bibitem[{Adams et~al.(2021)Adams, Colella, Graves, Johnson, Johansen, Keen,
  Ligocki, Martin, McCorquodale, Modiano, Schwartz, Sternberg, and
  Straalen}]{Chombo}
\bibinfo{author}{M.~Adams}, \bibinfo{author}{P.~Colella},
  \bibinfo{author}{D.~Graves}, \bibinfo{author}{J.~Johnson},
  \bibinfo{author}{H.~Johansen}, \bibinfo{author}{N.~Keen},
  \bibinfo{author}{T.~Ligocki}, \bibinfo{author}{D.~Martin},
  \bibinfo{author}{P.~McCorquodale}, \bibinfo{author}{D.~Modiano},
  \bibinfo{author}{P.~Schwartz}, \bibinfo{author}{T.~Sternberg},
  \bibinfo{author}{B.~V. Straalen}, \bibinfo{title}{Chombo software package for
  AMR applications: design document}, \bibinfo{type}{Technical Report},
  \bibinfo{year}{April 2021}.
\bibitem[{Childs et~al.(2012)Childs, Brugger, Whitlock, Meredith, Ahern,
  Pugmire, Biagas, Miller, Harrison, Weber, Krishnan, Fogal, Sanderson, Garth,
  Bethel, Camp, Rubel, Durant, Favre, and Navratil}]{visit}
\bibinfo{author}{H.~Childs}, \bibinfo{author}{E.~Brugger},
  \bibinfo{author}{B.~Whitlock}, \bibinfo{author}{J.~Meredith},
  \bibinfo{author}{S.~Ahern}, \bibinfo{author}{D.~Pugmire},
  \bibinfo{author}{K.~Biagas}, \bibinfo{author}{M.~C. Miller},
  \bibinfo{author}{C.~Harrison}, \bibinfo{author}{G.~H. Weber},
  \bibinfo{author}{H.~Krishnan}, \bibinfo{author}{T.~Fogal},
  \bibinfo{author}{A.~Sanderson}, \bibinfo{author}{C.~Garth},
  \bibinfo{author}{E.~W. Bethel}, \bibinfo{author}{D.~Camp},
  \bibinfo{author}{O.~Rubel}, \bibinfo{author}{M.~Durant},
  \bibinfo{author}{J.~M. Favre}, \bibinfo{author}{P.~Navratil},
  \bibinfo{title}{{VisIt: An End-User Tool For Visualizing and Analyzing Very
  Large Data}}, \bibinfo{year}{2012}.
\bibitem[{Bochkov and Gibou(2020)}]{gibou}
\bibinfo{author}{D.~Bochkov}, \bibinfo{author}{F.~Gibou},
\newblock \bibinfo{title}{Solving elliptic interface problems with jump
  conditions on cartesian grids},
\newblock \bibinfo{journal}{Journal of Computational Physics}
  \bibinfo{volume}{407} (\bibinfo{year}{2020}) \bibinfo{pages}{109269}.

\end{thebibliography}
